\documentclass[12pt]{amsart}

\usepackage[margin=1.1in]{geometry}

\usepackage{url}

\usepackage{amsfonts}
\usepackage{amsthm,amssymb}
\usepackage{mathtools}

\usepackage[hidelinks]{hyperref} 
\usepackage[capitalize, nameinlink, noabbrev]{cleveref}

\usepackage{graphicx}
\usepackage{subcaption}
\graphicspath{ {./figs/}{../figs/} }

\usepackage[dvipsnames]{xcolor}

\newcommand{\defref}[1]{Definition~\ref{def:#1}}
\newcommand{\figref}[1]{Figure~\ref{fig:#1}}

\newcommand{\thmref}[1]{Theorem~\ref{thm:#1}}
\newcommand{\defcondref}[2]{Condition~\ref{defcond:#1-#2} of \defref{#1}}

\usepackage{tikz}
\usepackage{tikz-cd}
\usetikzlibrary{graphs,arrows,decorations.pathmorphing,decorations.markings,fit,positioning,hobby,arrows.meta}
\usepackage[scr]{rsfso}

\tikzset{
    marrow/.style={decoration={markings,mark=at position 0.5 with {\arrow{#1}}}, postaction=decorate}
}
\tikzset{%
vcenter/.style = {baseline = {([yshift = -.5ex](0, 1.25))}},%
functor fill/.style = {fill = black!30!white, opacity = .6},%
functor draw/.style = {line width = 2pt, draw = black!50!white, opacity = .6}
}

\def\R{\mathbb{R}}

\DeclareMathOperator{\codis}{codis}
\DeclareMathOperator{\dis}{dis}
\DeclareMathOperator{\id}{id}

\DeclareMathOperator{\len}{len}

\renewcommand{\P}{\vec{P}} 
\def\PX{\P(X)}
\def\PY{\P(Y)}

\def\PZZ{\P_{zz}} 
\def\X{\vec{X}} 
\def\Y{\vec{Y}}
\def\Z{\vec{Z}}
\def\dF{\vec{F}} 
\def\dG{\vec{G}}
\def\dZ{d_{zz}} 
\def\lZ{l_{zz}}

\def\diam{\mathrm{Diam}}

\newcommand{\haus}{d_\mathrm{H}}
\newcommand{\dhaus}{\vec{d}_\mathrm{H}}
\newcommand{\GH}{d_\mathrm{GH}}
\newcommand{\dGH}{\vec{d}_\mathrm{GH}}
\newcommand{\ddis}{\vec{d}_{\textup{dis}}}
\newcommand{\Xrev}{\X^\ast}

\newcommand{\define}[1]{\textbf{#1}}

\theoremstyle{definition}
\newtheorem{theorem}{Theorem}
\newtheorem{corollary}{Corollary}[section]
\newtheorem{lemma}{Lemma}[section]
\newtheorem{proposition}{Proposition}[section]
\newtheorem{definition}{Definition}[section]
\newtheorem{remark}{Remark}[section]
\newtheorem{example}{Example}[section]

\usepackage{color}
\definecolor{darkorchid}{rgb}{0.6,0.196,0.8}

\newcommand\TODO[3]{\hbox to 0pt{\textcolor{#1}{$^\bullet$}}\marginpar{\footnotesize \textcolor{#1}{#2: #3}}}

\title{Gromov--Hausdorff Distance for Directed Spaces}

\author[L.Fajstrup]{Lisbeth Fajstrup}
\address{Aalborg University, 9220 Aalborg {\O}st, Denmark}
\email{fajstrup@math.aau.dk}

\author[B.T.Fasy]{Brittany Terese Fasy}
\address{Montana State University,
    Bozeman, MT 59715, USA}
\email{brittany.fasy@montana.edu}

\author[W.Li]{Wenwen Li}
\address{Hobart and William Smith Colleges, Geneva, NY 14456, USA}
\email{wli11uco@gmail.com}

\author[L.Mezrag]{Lydia Mezrag}
\address{Université de Montréal, Montréal, QC  6128, CA}
\email{lydia.mezrag@umontreal.ca}

\author[T.Rask]{Tatum Rask}
\address{Colorado State University, 
    Fort Collins, CO 80524, USA}
\email{tatum.rask@colostate.edu}

\author[F.Tombari]{Francesca Tombari}
\address{Max Planck Institute for Mathematics in the Sciences, Leipzig 04103, Germany}
\email{francesca.tombari@mis.mpg.de}

\author[\v{Z}.Urban\v{c}i\v{c}]{\v{Z}iva Urban\v{c}i\v{c}}
\address{Durham University,
    Durham, DH1 3LE, UK}
\email{ziva.urbancic@durham.ac.uk}

\begin{document}
\maketitle          

\begin{abstract}
The Gromov--Hausdorff distance measures the similarity between two metric spaces
by isometrically embedding them into an ambient metric space. 
We introduce an analogue of this distance for metric spaces endowed with directed structures. 
The directed Gromov--Hausdorff distance measures the distance between two extended metric spaces, where the new metric, defined on the same underlying space, is induced by the length of zigzag paths.
This distance is then computed by isometrically embedding the directed metric spaces into an ambient directed space equipped with the zigzag distance.
Analogously to the classical Gromov--Hausdorff distance, we also propose alternative formulations based on the distortion of d-maps and d-correspondences. However, unlike the classical case, these directed distances are not equivalent.
\end{abstract}

\section{Introduction}\label{sec:intro}

Directed algebraic topology emerged in the 1990s as a framework for modelling
non-reversible phenomena.
It has been developed in various forms, notably in homotopy theory~\cite{grandis1993cubical} and as a theoretical model for concurrency
\cite{lisbeth,gaucher2000homotopy}. 
More recently, the increasing use of networks in mathematical modeling has drawn renewed
attention to the study of directed structures in different areas, such as
applied topology~\cite{masulli2016,Chowdhury_PersistentPathHomology},
combinatorics~\cite{dochtermann2023homomorphism}, and machine learning
\cite{Reimann_2017,zia2025persistent}.
In this paper, a directed space refers to a 
topological space equipped with a distinguished set of continuous paths, called directed
paths, which satisfy specific axioms.
Examples of directed spaces include directed graphs and topological spaces
with partial orders.
As a result of this line of research, new distances, sensitive to
directionality, have been introduced, see for example~\cite{chowdhury2018}
and~\cite{zava2019}.
Taking directionality into account is crucial not only for distinguishing
networks, but also for establishing stability results and measuring the robustness of feature
detection in data analysis\cite{zhang2024}.

The Gromov--Hausdorff distance provides a way to compare
metric spaces by embedding them isometrically into a common ambient
space and measuring the extent to which they fail to be isometric.
It was introduced in the 1970s as an extension of the Hausdorff
distance~\cite{edwards1975structure,gromov1981groups,gromov1981structures}.
In fact, the Gromov--Hausdorff distance defines a metric on the space of compact
metric spaces up to isometry \cite{ivanov2015gromov}. Applications span
various fields, including the comparison of metric
graphs~\cite{agarwal2018computing,lee2011computing}
and shape analysis~\cite{bronstein2008numerical,memoli2004comparing}.
Inspired by the classical Gromov--Hausdorff distance, we define in this paper the \emph{directed Gromov--Hausdorff distance}, denoted~$\dGH$, on the set of
metric spaces $(X,d^X)$ endowed with a directed structure~$\PX$. 
This construction involves introducing a new notion of isometry on d-spaces, called
\emph{d-isometry}, and computing a variant of the Hausdorff distance with
respect to an extended metric called the \emph{zigzag distance}, denoted
$\dZ$.
The latter is constructed as an induced extended metric by considering the
length of admissible directed paths connecting points in $X$. 

The construction of the zigzag distance is closely related to the notion
of \emph{length}, a well-studied topic in metric geometry. A
\emph{length structure} on a topological space $X$ consists of a collection $\mathcal{C}$ of
continuous paths and a positive map~$\len : \mathcal{C} \rightarrow [0,
\infty]$ satisfying a set of axioms that differ depending on the source.
We refer to \cite[\textsection~2.1]{burago2001course} for an example of such
axioms.
In our work, we restrict ourselves to metric spaces and consider only rectifiable paths as
admissible, ensuring that the length of a path and, consequently, the zigzag distance
are well-defined and~finite.

While investigating some properties of the directed Gromov--Hausdorff distance, we recall equivalent formulations of the classical Gromov--Hausdorff distance based on distortion \cite{burago2001course,kalton1999}.
By extending these constructions to the directed setting, we introduce new distances between directed spaces: the \emph{directed distortion distance} and the \emph{d-correspondence distortion distance}. 
Interestingly, while these distances are equivalent to the Gromov--Hausdorff distance for (undirected) metric spaces, this equivalence does not hold for directed metric spaces. 

The structure of the paper is as follows.
In \cref{sec:prelim}, we provide the necessary background information on
directed topology and the classical Gromov--Hausdorff distance. 
Then, in \cref{sec:directed-GH}, we introduce directed versions of the Gromov--Hausdorff
distance. Specifically, we
start in \cref{sub:zigzag-distance} by describing
how a metric space~$(X, d^X)$
endowed with a directed structure $\PX$ induces an extended metric~$\dZ^X$ on
the directed space~$\X$ and explore the topology induced by this metric. 
Next, in \cref{sub:directed-GH-distance}, we define the directed
Gromov--Hausdorff distance between two directed metric spaces and explore some
of its properties. In \cref{sub:distortion-distance}, we define a directed
analogue of distorsion-based distances between metric spaces. 

\section{Preliminaries}\label{sec:prelim}

In this section, we recall definitions and results on directed spaces and
the classical Gromov--Hausdorff distance.
For a general introduction to these topics, we refer to \cite{lisbeth} and \cite{Grandis_2009} for the
former, and to \cite[\textsection~7.3]{burago2001course} for the latter.

\subsection{Directed Spaces}

Denote by $I$ the unit interval $[0,1]$ with the standard topology and consider another topological space $X$.
A continuous map $\gamma\colon  I\to X$ is called a \define{path}
in~$X$; the points~$\gamma(0)$ and~$\gamma(1)$ are called the \define{source}
and \define{target} of the path,
respectively.
Specifying the source and the target of a path endows it with a notion of
direction, which allows us to say $\gamma$ is a path from $s$ to $t$.
Given two paths, $\gamma$ and $\gamma'$, such that $\gamma$ has source~$x$ and target~$x'$, and $\gamma'$ has
source $x'$ and target~$x''$, their \define{concatenation}, $\gamma\gamma'\colon I \to X$, defined as
\begin{equation*}
    \gamma\gamma'(t) :=
        \begin{cases}
            \gamma(2t)      & 0\le t\le \frac{1}{2}, \\
            \gamma'(2t-1)   & \frac{1}{2}\le t\le 1,
        \end{cases}
\end{equation*}
is a path with source $x$ and target~$x''$.
A path $\gamma \colon I \rightarrow X$ is a \define{reparametrization} of $\mu \colon I\to
X$ if there is a nondecreasing, continuous, surjective map $h \colon I\to I$ such that
$\gamma=\mu\circ h$.

\begin{remark}[The Unit Interval]\label{rem:paths-defined-on-unit-interval}
    Our paths are defined on the unit interval following the tradition in
    homotopy theory. As a consequence, the concatenation is not associative on
    the nose, but up to reparametrization,
    which is all we need.
\end{remark}

\begin{remark}[Taking Subpaths]
    Given a path~$\gamma \colon I\to X$, the
    restriction $\mu:=\gamma|_{[a,b]}$ of the path to a
    subinterval~$[a,b]\subseteq I=[0,1]$ is a path if and only if
    $[a,b]=I$, as, by definition, the domain of a
    path must be $I$. However,
    by precomposing~$\mu$ with the linear rescaling function~$r \colon I \to [a,b]$
    defined by~$r(t)=a+t(b-a)$, the map~$\mu\circ r$ is a path such that~$\mu
    \circ r(0)=\mu(a)$,~$\mu\circ r(1)=\mu(b)$ and $\mu \circ
    r(I)=\mu([a,b])$.
    We refer to $\mu\circ r$ as a \define{subpath} of~$\gamma$.
\end{remark}

\begin{definition}\label{def:dspace}
    A \define{directed space} or \define{d-space} is a pair $(X, \PX)$, where
    $X$ is a topological space and $\PX$ is a collection of paths on $X$ such
    that
    \begin{enumerate}
        \item every constant path belongs to
            $\PX$,\label{defcond:dspace-constant}
        \item $\PX$ is closed under
            reparameterization,\label{defcond:dspace-reparam}
        \item $\PX$ is closed under taking
            subpaths,\label{defcond:dspace-subpath}
        \item if $\gamma,\gamma'$ are two paths in $\PX$ such that $\gamma(1)=\gamma'(0)$, then
            the concatenation~$\gamma \gamma'$ is
            also in~$\PX$.\label{defcond:dspace-concat}
    \end{enumerate}
    We denote the pair $(X, \PX)$ as $\X$ whenever the set of
    d-paths is clear from context.
    Furthermore, we refer to $\PX$ as the set of \define{d-paths} of $\X$ or the
    \define{d-structure} on $X$, and we
    denote by $\P(x,x')$ the subset of $\PX$ containing those d-paths with
    $x$ as the source and $x'$ as the target.
\end{definition}
Combined, Conditions~\ref{defcond:dspace-reparam}
and~\ref{defcond:dspace-subpath} imply that
$\PX$ contains all
reparameterizations of all subpaths of $\gamma$.

Given two d-spaces $\X$ and $\Y$, a \define{d-map}, or directed map, $\dF\colon \X\to \Y$ is a
continuous map on the underlying spaces~$F\colon X\to Y$ such that, for every
path~$\gamma$ in $\PX$, the
composition~$F\circ\gamma$ is in $\PY$.

A topological space can be endowed with many d-structures.
For example, the d-space $(X, \PX)$ is called the \define{discrete} d-space on $X$
if~$\PX$ contains only the constant paths, $\PX=\{ c_x \}_{x \in X}$, and the
\define{trivial} or indiscrete d-structure on~$X$ if it contains all
continuous maps from $I$ to $X$,
$\PX=C(I, X)$.

\begin{remark}\label{rem:dTop}
The category \define{dTop} is the category with d-spaces as objects and d-maps
as morphisms.
This category is complete and cocomplete; see \cite[Prop.~4.5]{lisbeth}.
The discrete d-structure is a left adjoint to the forgetful
functor from \define{dTop} to \define{Top} and the trivial d-structure is
a right adjoint to the forgetful functor. The analogy to the discrete topology
is that the discrete topology is a left adjoint to the forgetful functor from \define{Top} to
\define{Set}. Similarly, the trivial or indiscrete  topology is a right~adjoint.
\end{remark}

If $\{(X,\vec P_j(X))\}_{j\in J}$ is a collection of d-spaces, then $(X, \bigcap_j \vec P_j(X))$ is also a d-space.
On the other hand, $(X, \bigcup_{j}\vec P_j(X))$ need not be a d-space, as the
set of paths might not be closed under concatenation
(\defcondref{dspace}{concat}).
However, we talk about the d-structure generated by $\bigcup_{j}\vec P_j (X)$.
Given a set~\mbox{$A \subset \PX$}, we say that $A$ \define{generates} $\PX$ if $\PX$
is the smallest of all d-structures containing~$A$;
see~\cite[Def.~3.8]{fajstrup2017hierarchy}.
In this case, we write~$\PX=\langle A\rangle$.

The pair $(Y, \PY)$ is a \define{d-subspace} of the directed space~$\X=(X, \PX)$ if $Y\subset X$
and~$\PY=\{\gamma\in\PX | \gamma(I)\subset Y\} $.
In this way, every subset of~$X$ has a d-structure induced  by the one on $\X$.
Moreover, there is a
canonical d-map~$\vec{\iota}\colon \Y\to \X$ induced by the inclusion
$\iota\colon Y\subset X$.

For two d-spaces $\X$ and $\Y$, a bijective d-map $\dF\colon \X\to \Y$ whose inverse is a d-map is called
\define{d-invertible}.
For example, if the map $F \colon X \to Y$ is a
homeomorphism and if $\dF$ is bijective on d-paths, then $\vec F$ is d-invertible.
However,
not all bijective d-maps are d-invertible, even if they are homeomorphisms on
the underlying spaces:
the identity d-map from $(I, \P_1(I))$ to~$(I, \P_2(I))$,
where $\P_1(I)$ is the discrete structure and $\P_2(I)$
is the trivial structure,
is not d-invertible.

\subsubsection{Path Length}

Consider a metric space $(X, d)$ and endow $X$ with the topology induced by the
metric.
For a path $\gamma \colon I \to X$, its \define{length} is defined as
\begin{equation}\label{eqn:length}
    \len(\gamma)=\sup \sum_{i=1}^N d(\gamma(t_{i-1}), \gamma(t_i)),
\end{equation}
where the supremum ranges over all $N\in\mathbb{N}$ and all sequences $0\le
t_0<t_1<\cdots<t_N\le 1$; see~\cite{sullivan2008curves}.
We say that $\gamma$ is a \define{rectifiable}
path if $\len(\gamma)<\infty$.

\begin{definition}\label{rectifiabledspace}
    A \define{rectifiable d-space} is a d-space $(X, \PX)$, where $X$ is a
    metric space and every d-path in $\PX$ is rectifiable.
\end{definition}

\begin{remark}\label{rem:rec-pseudo}
    Rectifiability is defined for paths in pseudo-metric spaces and extended
    metric spaces\footnote{While the definitions of pseudo-metric and extended
    metric are not common across fields, we adopt the terminology used
    in computational geometry~\cite{conci2018}. Thus, a
    pseudo-metric follows all metric properties, except separability;
    that is, $d(x, x')=0$ does not necessarily imply that $x=x'$.
    Furthermore, an extended metric follows all metric properties, except
    finiteness.} in the same way.
    These generalized distances still induce a topology and we extend
    \cref{rectifiabledspace} to such spaces.
\end{remark}

Note that every constant path is rectifiable (as it has length zero)
and the concatenation of a finite
sequence of rectifiable paths is still rectifiable. 
A subpath of a rectifiable path is also rectifiable.
Moreover, if~$\gamma$ is rectifiable and~$h\colon I\to I$ is nondecreasing,
then $\gamma \circ h$ is rectifiable.
This is seen by observing that the collection of all sequences~$\{t_0<t_1<\dots<t_N\}$
contains the set of all sequences~$\{h(t_0)\le h(t_1)\le \dots\le h(t_N)\}$,
thus, looking at the supremum in \cref{eqn:length}, we find
\begin{equation}
    \sup \sum_{i=1}^N d(\gamma(h(t_{i-1})), \gamma(h(t_i)))
    \le \sup \sum_{i=1}^N d(\gamma(t_{i-1}), \gamma(t_i))
    <\infty.
\end{equation}
These observations guarantee that for a subset~$A\subseteq C(I, X)$ of rectifiable
paths, the d-structure $\langle A \rangle$ is rectifiable.
In addition, a d-subspace of a rectifiable d-space is~rectifiable.

\begin{example}[Partially Ordered Space]\label{ex:partial}
    Let $(X, \le)$ be a topological space with a partial order $\le$,
    such that $\{(x,x')\mid x\le x'\}\subset X\times X$ is closed in the
    product topology; see~\cite{lisbeth},
    and let $\PX$ be the subset of nondecreasing paths, i.e.,
    those paths $\gamma$ such that $\gamma(s)\le \gamma(t)$, for every $s\le t
    \in [0,1]$.
    Then, $(X, \PX)$ is a d-space.
    For example, $\R^n$, endowed with the product order, and all the
    nondecreasing paths is a d-space denoted by $\vec{\R}^n$.
    If, in addition,~$X$ is a metric space, the rectifiable
    nondecreasing paths define a rectifiable d-space. In particular, this holds
    for $\vec{\R}^n$.
\end{example}

\begin{example}[Reversed Path Space]\label{ex:reverse}
    Let $\X=(X, \PX)$ be a (rectifiable) d-space.
    Define the reversed
    (rectifiable) d-space $\Xrev:=(X,\PX^\ast)$,
    where $\PX^\ast$ is the set of
    d-paths obtained by reversing the direction of the paths in~$\PX$,
    i.e., $\gamma$ is in~$\PX$ if and only if $\gamma^\ast$, defined by
    $\gamma^\ast(t)=\gamma(1-t)$, is in~$\PX^\ast$.
\end{example}

\subsubsection{Some Operations on Directed Spaces}\label{subsec:operations}

Let $\X = (X, \PX)$ and $\Y = (Y, \PY)$ be two d-spaces. We
outline below a few constructions on d-spaces.

\paragraph{Cartesian Product} 
The cartesian product $\X \times \Y$ is defined as the topological
space $X \times Y$ together with the directed structure given by the product
$\P(X \times Y) := \PX \times \PY$, i.e., a path in $X\times Y$ is in
$\P(X\times Y)$ if and only if its components are, respectively, in $\PX$ and
$\PY$; see \cite[Section 1.4.1]{Grandis_2009}.

\paragraph{Disjoint Union}
The disjoint union $\X \sqcup \Y$ is
the topological space~$X \sqcup Y$ together with the directed
structure given by the coproduct $\P(X \sqcup Y) := \PX \sqcup
\PY$, i.e., a path in~$X\times Y$ is in~$\P(X\times Y)$ if and only if
either of its components is in~$\PX$ or in~$\PY$; see \cite[Section 1.4.1]{Grandis_2009}.
If $(X,d^X)$ and $(Y,d^Y)$ are \mbox{(pseudo-)metric} spaces, there are several ways to
induce distances (and metrics) on $X \sqcup Y$. For instance, consider the
extended (pseudo-)metric space $X \sqcup Y$ with distances:
\begin{equation*}
    d(a,b)=
    \begin{cases}
        d^X(a,b),   &\text{if } a,b\in X,\\
        d^Y(a,b),   &\text{if } a,b\in Y,\\
        \infty,     &\text{otherwise.}
    \end{cases}
\end{equation*}
If $\X$ and
$\Y$ are rectifiable d-spaces, then $\P(X \sqcup Y) := \PX \sqcup \PY$ gives a
rectifiable d-structure on~$X\sqcup Y$.

\paragraph{Quotient}
Let $\X$ be a directed space and let $\sim$ be an equivalence relation on the
underlying space~$X$.
The quotient d-space $\X/\sim$, with the quotient topology, is the d-space
whose d-paths,~$\P(X/\sim)$, form the smallest d-structure
such that the quotient map $\vec{\pi} \colon \X \to \X/\sim$ is a d-map. In
other words,~$\P(X/\sim)$ is generated by the set $ \vec{\pi}(\P) := \{ \vec{\pi} \circ \gamma \; | \; \gamma \in
\PX\}\subset C(I, \X/\sim)$.

If $\X$ is a rectifiable d-space and $\sim$ an equivalence relation, then we
define a pseudo-metric $\bar{d}$ on the quotient $X/\sim$ as
\begin{equation}
    \bar{d}([x],[y])
    := \inf\left\{\sum_{i = 1}^k d(x_i,y_i)
        \mid x_1 \in [x], y_k \in [y], x_{i+1} \sim y_i, k \in \mathbb{N}\right\}
\end{equation}
as found in \cite[Def. 3.1.12]{burago2001course}.

Let $\P(X/\sim)$ be the d-structure as generated above. To see that this defines
a rectifiable d-structure, it suffices to see that if $\gamma:I\to X$ is
rectifiable, then~$\vec{\pi}\circ \gamma$ is rectifiable. This follows from the fact
that for all $x,y \in X$, by definition,~$\bar{d}([x],[y])\leq d(x,y)$.  Thus,
for any sequence $\{t_0<\dots <t_N\}$,
we have~$\sum_{i=1}^N\bar{d}([\vec{\pi}\circ\gamma(t_{i-1})],
[\vec{\pi}\circ\gamma(t_i)]) \leq \sum_{i=1}^N d(\gamma(t_{i-1}), \gamma(t_i))$.

\begin{example}[Directed Graphs] \label{ex:dirgraph}
    Let~$G = (V, E, s, t)$ be a directed graph, where~$V$ is the set of
    vertices,~$E \subseteq V \times V$
    is the set of edges such that~$(v, v) \notin E$ for any~$v \in V$, and the maps~$s,t
    \colon E \to V$ assign the source (map~$s$) and target vertex (map~$t$) to
    each edge. We show here how to interpret this directed graph as a d-space. The
    requirements that a path be a continuous map and that the set of d-paths be
    closed under taking subpaths forces
    one to include edges as part of the underlying topological
    space. We leverage the fact that an undirected
    graph is a topological space, through the CW-complex construction,
    which we modify to impose a directed structure on it as well. Therefore,
    construct a topological space~$X_G$ by specifying its strata:
    \begin{itemize}
        \item $X_G^0 = V$,
        \item $X_G^1$ is formed by starting with~$X_G^0$ and attaching a one-cell for each~$e \in E$ via a
            gluing map~$\varphi_e \colon \partial \vec{I} \to X_G^0$
            with~$\varphi_e(0) = s(e)$ and~$\varphi_e(1) = t(e)$,
            where~$\vec{I}$ is the unit interval with the direction inherited
            from the usual partial order.
    \end{itemize}
    Consider the topological space defined by the disjoint union of
    $X_{G}^{0}$ and the one-cells in $X_{G}^{1}$.  We equip this space with the directed
    structure generated by constant paths and paths $\gamma_e$ such that
    $\gamma_e(0) = s(e)$ and $\gamma_e(1)=t(e)$, where~$e$ is an edge in $E$.
    Since~$X_G$ is the quotient space obtained from this space by imposing the
    attaching relation defined by the gluing maps $\varphi_e$, it inherits a
    d-structure from the d-structure described above.
\end{example}
\begin{example}[Directed Hollow Hypercubes]
    Let~$Q_n = \partial I^n \subset \mathbb{R}^n$ be the boundary of
    the~$n$-cube,
    and let~$\vec{\R}^n$ be~$\R^n$ with d-structure induced
    on it by the product order, as defined in \cref{ex:partial}.
    The directed hollow $n$-cube~$\vec{Q}_n$ is the d-space~$(Q_n,
    \P(Q_n))$ with the subspace directed structure inherited from~$\vec{\R}^n$.

    If we consider the rectifiable d-structure on ~$\vec{\R}^n$, the hollow
    hypercube inherits a rectifiable d-structure.
\end{example}

From now on, the only d-spaces we consider are going to be rectifiable, and so, we
just refer to them as d-spaces. The topology is the one induced by the metric,
pseudo-metric, or extended metric.
The reason for this choice is that we are interested in defining metrics on
d-spaces and, for this, having a notion of length of a path is paramount.

\subsection{Gromov--Hausdorff Distance}

Given a metric space $(M, d)$ and two nonempty subspaces $A, B\subset M$, the
Hausdorff distance between them is defined as
\[
    \haus(A,B)= \max\left\{\sup _{a\in A} d(a,B), \sup _{b\in B} d(A,b)\right\},
\]
where $ d(a,B)=\inf_{b\in B} d(a,b)$ and $d(A,b)=\inf_{a\in A}d(a, b)$.

Now, given two metric spaces $(X, d^X)$ and $(Y, d^Y)$, the
\define{Gromov--Hausdorff distance} between them is defined as
\[
    \GH(X, Y)=\inf_{f,g} \haus(f(X), g(Y)),
\]
where the infimum ranges over all metric spaces $(Z, d^Z)$ and over all
isometries~$f\colon X\hookrightarrow
Z$ and $g\colon Y\hookrightarrow Z$.

A relation between two sets $X$ and $Y$ is a
subset of $X\times Y$.  We call a relation~$\mathcal{R}$ a \define{correspondence}
if, for every $x$ in $X$ there exists $y$ in $Y$
such that~$(x, y)$ is in~$\mathcal{R}$, and, for every $y$ in $Y$,
there exists $x$ in $X$ such that~$(x, y)$ is in $\mathcal{R}$.
If~$(X, d^X)$ and $(Y, d^Y)$ are metric spaces and $\mathcal{R}\subset X\times Y$
is a non-empty relation, we define the \define{distortion} of $\mathcal{R}$ as
\[
    \dis(\mathcal{R})
    =\sup_{(x, y),(x', y')\in \mathcal{R}}
        \left\lvert d^X(x,x')-d^Y(y, y')\right\rvert.
\]

Observe that every function $f \colon X\to Y$ is associated with a relation,
$\mathcal{R}_f$, whose elements are pairs $(x, f(x))$.
Note that such relation $\mathcal{R}_f$ is a correspondence if and only if $f$
is surjective.
We denote the distortion $\dis(\mathcal{R}_f)$ by $\dis(f)$ and simply call it
the distortion of $f$.
Furthermore, given two functions~$f\colon X\to Y$ and $g\colon Y\to X$, where
$(X, d^X)$ and $(Y, d^Y)$ are metric spaces, the \define{codistortion} of the
pair $(f, g)$ is defined as
\[
    \codis(f, g)
    =\sup_{x\in X, \; y\in Y} \left\lvert d^X(x, g(y))-d^Y(f(x), y) \right\rvert.
\]
From~\cite[Theorem 7.3.25]{burago2001course}, we know that
\begin{equation} \label{eqn:kalton_bis}
    \GH(X, Y)
    =\frac{1}{2}\inf_{\mathcal{R}} \dis(\mathcal{R}),
\end{equation}
where $\mathcal{R}$ ranges over all the correspondences between $X$ and $Y$.
It was also observed in \cite{kalton1999} that, for $X$ and $Y$ bounded,
\begin{equation} \label{eqn:kalton}
    \GH(X, Y)
    =\frac{1}{2}\inf_{f, g} \max\{\dis(f), \dis(g), \codis(f, g)\},
\end{equation}
where $f\colon X\to Y$ and $g\colon Y\to X$ are (not necessarily
continuous) functions.

The following proposition collects some well-known results about the
Gromov--Hausdorff distance, which are not hard to show.
Given a metric space $(X, d)$, its diameter $\diam(X, d)$
is defined as ${\sup_{x,x'\in X} d(x,x')}$.

\begin{proposition}[Properties of the Gromov--Hausdorff Distance \cite{burago2001course}]\label{prop:GH_properties}
    Let~$(X,d^X)$ and $(Y, d^Y)$ be metric spaces. Then, the following holds:
    \begin{enumerate}
        \item $d_{GH}(X,Y) \leq \frac{1}{2} \sup \lbrace \diam(X,d^X), \diam(Y, d^Y)\rbrace$.
        \item $d_{GH}(X,Y) < \infty$, if $X$ and $Y$ are bounded.
        \item $d_{GH}(X,Y) = \frac{1}{2}\diam(Y, d^Y)$, if $X=\{x_0\}$.
        \item $d_{GH}(X,Y)\ge \frac{1}{2} |\diam(X,d^X)-\diam(Y, d^Y)| $,  if $\diam(X,d^X) < \infty $.
        \item $d_{GH}(X,Y) \geq \frac{1}{2}\inf\{\dis (f)\mid f\colon X\to Y\}$.\label{GH_prop_4}
    \end{enumerate}
\end{proposition}

\section{Directed Gromov--Hausdorff Distance}\label{sec:directed-GH}

In this section, we begin by introducing the zigzag distance on a metric space
with a directed structure and showing that it is, indeed, an extended metric.
This distance is similar in spirit to the path (or intrinsic) metric of
a metric space, as presented in, e.g.,~\cite{gromov1999metric},
but it is a new application of these ideas to directed spaces.
We then use this new distance to define an analogue of the Gromov--Hausdorff distance for d-spaces.

\subsection{Zigzag Distance}\label{sub:zigzag-distance}

A \define{zigzag path} between $x, x' \in \X$ is a sequence of d-paths~$(
\gamma_i )_{i = 1}^m$ such that, for each $i \in\{1,2,\ldots, m\}$, we have~$\gamma_i \in
\P(p_{i-1}, p_i) \cup \P(p_i, p_{i-1})$, where $x=p_0$ and~$x'=p_m$; see \figref{zigzag_paths} for an example.
Denote the set of all zigzag paths between $x$ and $x'$ as $\PZZ(x,
x')$, and
the set of all zigzag paths in $\X$
as~$\PZZ(X)$.
Assume~$X$ is a metric space, and recall that we assume the d-paths are rectifiable.
We define the length of a zigzag path~$\gamma = ( \gamma_i)_{i = 1}^m$~as
\begin{equation}
    \lZ(\gamma) = \sum_{i = 1}^{m} \len(\gamma_i).
\end{equation}
A d-space $\X$ is \define{zigzag connected} if 
$\PZZ(x, x') \neq \varnothing$ for all $x, x' \in X$.
A zigzag connected d-space $\X$ is also path connected because $\PZZ(X)\subset
C([0,1],X)$.
The converse, however, may not be true.
As an example, consider a path connected space with $\PX$ containing
only the constant paths.

\begin{figure}[ht]\label{fig:zigzag_paths}
        \centering
        \begin{subfigure}{.25\textwidth}
            \centering

\tikzset{every picture/.style={line width=0.75pt}} 

\begin{tikzpicture}[x=0.75pt,y=0.75pt,yscale=-1,xscale=1]

\draw  [fill={rgb, 255:red, 155; green, 155; blue, 155 }  ,fill opacity=0.5 ] (281.6,21) -- (376,21) -- (376,115.4) -- (281.6,115.4) -- cycle ;
\draw [color={rgb, 255:red, 126; green, 211; blue, 33 }  ,draw opacity=1 ][line width=1.5]    (281.6,115.4) -- (373.88,23.12) ;
\draw [shift={(376,21)}, rotate = 135] [color={rgb, 255:red, 126; green, 211; blue, 33 }  ,draw opacity=1 ][line width=1.5]    (14.21,-4.28) .. controls (9.04,-1.82) and (4.3,-0.39) .. (0,0) .. controls (4.3,0.39) and (9.04,1.82) .. (14.21,4.28)   ;
\draw [color={rgb, 255:red, 126; green, 211; blue, 33 }  ,draw opacity=1 ]   (281.6,115.4) -- (376,21) ;
\draw [shift={(376,21)}, rotate = 315] [color={rgb, 255:red, 126; green, 211; blue, 33 }  ,draw opacity=1 ][fill={rgb, 255:red, 126; green, 211; blue, 33 }  ,fill opacity=1 ][line width=0.75]      (0, 0) circle [x radius= 3.35, y radius= 3.35]   ;
\draw [shift={(281.6,115.4)}, rotate = 315] [color={rgb, 255:red, 126; green, 211; blue, 33 }  ,draw opacity=1 ][fill={rgb, 255:red, 126; green, 211; blue, 33 }  ,fill opacity=1 ][line width=0.75]      (0, 0) circle [x radius= 3.35, y radius= 3.35]   ;
\end{tikzpicture}
            \caption{A d-path.}
            \label{subfig:d-path}
        \end{subfigure}
        \hfill
        \begin{subfigure}{.25\textwidth}
            \centering

\tikzset{every picture/.style={line width=0.75pt}} 

\begin{tikzpicture}[x=0.75pt,y=0.75pt,yscale=-1,xscale=1]

\draw  [fill={rgb, 255:red, 155; green, 155; blue, 155 }  ,fill opacity=0.5 ] (282.6,18.2) -- (377,18.2) -- (377,112.6) -- (282.6,112.6) -- cycle ;
\draw [color={rgb, 255:red, 126; green, 211; blue, 33 }  ,draw opacity=1 ][line width=1.5]    (377,112.6) -- (377,68.2) ;
\draw [shift={(377,65.2)}, rotate = 90] [color={rgb, 255:red, 126; green, 211; blue, 33 }  ,draw opacity=1 ][line width=1.5]    (14.21,-4.28) .. controls (9.04,-1.82) and (4.3,-0.39) .. (0,0) .. controls (4.3,0.39) and (9.04,1.82) .. (14.21,4.28)   ;
\draw [color={rgb, 255:red, 126; green, 211; blue, 33 }  ,draw opacity=1 ][line width=1.5]    (374,65.21) -- (329.8,65.4) ;
\draw [shift={(377,65.2)}, rotate = 179.76] [color={rgb, 255:red, 126; green, 211; blue, 33 }  ,draw opacity=1 ][line width=1.5]    (14.21,-4.28) .. controls (9.04,-1.82) and (4.3,-0.39) .. (0,0) .. controls (4.3,0.39) and (9.04,1.82) .. (14.21,4.28)   ;
\draw [color={rgb, 255:red, 126; green, 211; blue, 33 }  ,draw opacity=1 ][line width=1.5]    (329.8,65.4) -- (329.8,21) ;
\draw [shift={(329.8,18)}, rotate = 90] [color={rgb, 255:red, 126; green, 211; blue, 33 }  ,draw opacity=1 ][line width=1.5]    (14.21,-4.28) .. controls (9.04,-1.82) and (4.3,-0.39) .. (0,0) .. controls (4.3,0.39) and (9.04,1.82) .. (14.21,4.28)   ;
\draw [color={rgb, 255:red, 126; green, 211; blue, 33 }  ,draw opacity=1 ][line width=1.5]    (326.8,18.01) -- (282.6,18.2) ;
\draw [shift={(329.8,18)}, rotate = 179.76] [color={rgb, 255:red, 126; green, 211; blue, 33 }  ,draw opacity=1 ][line width=1.5]    (14.21,-4.28) .. controls (9.04,-1.82) and (4.3,-0.39) .. (0,0) .. controls (4.3,0.39) and (9.04,1.82) .. (14.21,4.28)   ;
\draw [color={rgb, 255:red, 126; green, 211; blue, 33 }  ,draw opacity=1 ]   (329.8,18) -- (282.6,18.2) ;
\draw [shift={(282.6,18.2)}, rotate = 179.76] [color={rgb, 255:red, 126; green, 211; blue, 33 }  ,draw opacity=1 ][fill={rgb, 255:red, 126; green, 211; blue, 33 }  ,fill opacity=1 ][line width=0.75]      (0, 0) circle [x radius= 3.35, y radius= 3.35]   ;
\draw [color={rgb, 255:red, 126; green, 211; blue, 33 }  ,draw opacity=1 ][line width=0.75]    (377,112.6) -- (377,65.2) ;
\draw [shift={(377,112.6)}, rotate = 270] [color={rgb, 255:red, 126; green, 211; blue, 33 }  ,draw opacity=1 ][fill={rgb, 255:red, 126; green, 211; blue, 33 }  ,fill opacity=1 ][line width=0.75]      (0, 0) circle [x radius= 3.35, y radius= 3.35]   ;
\end{tikzpicture}
            \caption{A zigzag path.}
            \label{subfig:zigzag_path}
        \end{subfigure}
        \hfill
        \begin{subfigure}{.25\textwidth}
            \centering

\tikzset{every picture/.style={line width=0.75pt}} 

\begin{tikzpicture}[x=0.75pt,y=0.75pt,yscale=-1,xscale=1]

\draw  [fill={rgb, 255:red, 155; green, 155; blue, 155 }  ,fill opacity=0.5 ] (283.6,22) -- (378,22) -- (378,116.4) -- (283.6,116.4) -- cycle ;
\draw [color={rgb, 255:red, 208; green, 2; blue, 27 }  ,draw opacity=1 ]   (378,116.4) -- (283.6,22) ;
\draw [shift={(283.6,22)}, rotate = 225] [color={rgb, 255:red, 208; green, 2; blue, 27 }  ,draw opacity=1 ][fill={rgb, 255:red, 208; green, 2; blue, 27 }  ,fill opacity=1 ][line width=0.75]      (0, 0) circle [x radius= 3.35, y radius= 3.35]   ;
\draw [shift={(378,116.4)}, rotate = 225] [color={rgb, 255:red, 208; green, 2; blue, 27 }  ,draw opacity=1 ][fill={rgb, 255:red, 208; green, 2; blue, 27 }  ,fill opacity=1 ][line width=0.75]      (0, 0) circle [x radius= 3.35, y radius= 3.35]   ;
\draw [color={rgb, 255:red, 208; green, 2; blue, 27 }  ,draw opacity=1 ][line width=1.5]    (378,116.4) -- (342.3,80.7) -- (285.72,24.12) ;
\draw [shift={(283.6,22)}, rotate = 45] [color={rgb, 255:red, 208; green, 2; blue, 27 }  ,draw opacity=1 ][line width=1.5]    (14.21,-4.28) .. controls (9.04,-1.82) and (4.3,-0.39) .. (0,0) .. controls (4.3,0.39) and (9.04,1.82) .. (14.21,4.28)   ;
\end{tikzpicture}
            \caption{\emph{Not} a zigzag path.}
            \label{subfig:zigzag_nonexample}
        \end{subfigure}
        \caption{Let $\vec{I}^2$ be the d-space where $\P(I^2)$ is given by the
        product order on $I^2$. The single green arrow (left) is a d-path from
        $(0,0)$ to $(1,1)$. The multiple green arrows (middle) give a zigzag
        path between $(1,0)$ and $(0,1)$, whereas the red arrow (right) is \emph{not} a zigzag path.}
        \label{fig:zigzag_paths}
    \end{figure}
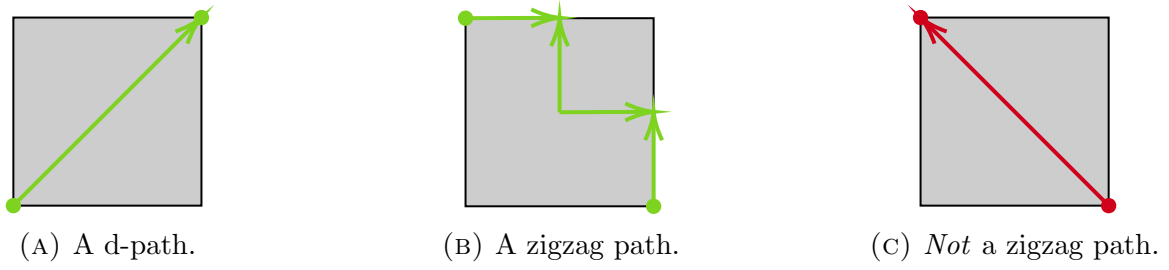

\begin{definition} \label{def: zigzag_length}
    Let $(X, d)$ be a metric space and $\X=(X,\PX) $ a d-space.
    For every $x,x'$ in $X$, define the \define{zigzag distance} induced by $d$ and $\PX$ on~$\X$ as
    \[
        \dZ(x, x')
            = \inf_{\gamma \in \PZZ(x, x')} \lZ(\gamma).
    \]
\end{definition}

\begin{remark}\label{rem:samedzz}
    Observe that different metric spaces with the same directed structure may define the same zigzag distance.
    This happens, for example, when considering the two metric spaces
    $(X, L^p)$ and $(X, L^q)$, with $p\neq q$ and $X=\{(x,1)\in \R^2\mid 0\le x\le 1 )\}\cup \{(1,x)\in \R^2\mid 0\le x\le 1 )\}$,
    both with the collection of paths $\vec P(X)$ given by all monotone
    nondecreasing paths.  On the other hand, it is possible to consider
    different d-structures on a metric space $(X, d)$ so that the associated
    zigzag distances are equal.  For example, consider $(X, \vec P(X))$ and its
    reversed path space $(X,\vec P(X)^*)$.
\end{remark}

\begin{lemma}\label{lem_zigzag_bigger_d}
    Let $\X = (X, d)$ be a metric space
    and $(\X, \dZ)$ be a d-space together with the zigzag distance $\dZ$ induced by $d$ and $\PX$.
    Then, for every~$x,x'$ in $X$, we have~$\dZ(x, x')\ge d(x, x')$.
\end{lemma}
\begin{proof}
    We first note that if $\PZZ(x, x')= \varnothing$, then $\dZ(x, x')=\infty$,
    which would mean that the
    statement holds.
    Otherwise, if $\PZZ(x, x')\neq\varnothing $, then
    for every zigzag path $\gamma = (\gamma_i)_{i=1}^m \in \PZZ(x,x')$ with
    $\gamma_1(0) = x$ and $\gamma_m(1)=x'$, we have~that
    \begin{equation}
        d(x,x') \leq \sum_{i=1}^m d(\gamma_i(0),\gamma_i(1)) \leq \sum_{i=1}^m
        \len(\gamma_i) \leq \lZ(\gamma).
    \end{equation}
    The first inequality follows from the triangle inequality.
    Taking the infimum over all zigzag paths $\gamma$, we conclude that $d(x,x') \leq \dZ(x, x')$.
\end{proof}

\begin{proposition}\label{prop_zigzag_metric}
    Let $(X,d)$ be a metric space together with the directed structure~$\vec
    P(x)$. Then, the zigzag distance, $\dZ$, on $\X$ is an extended metric.
\end{proposition}

\begin{proof}
    Because the infimum over the empty set is defined to be $\infty$,
    we know that if~$\PZZ(x,x')=\varnothing$, then $\dZ(x,x')=\infty$.
    If $\dZ(x, x')=0$ then, by \cref{lem_zigzag_bigger_d},
    we know that~$d(x, x')=0$.
    Because~$d$ is a metric, we conclude that $x=x'$.
    Now, observe that $\PZZ(x, x')=\PZZ(x', x)$,
    hence~$\dZ(x, x')=\dZ(x', x)$, implying that~$\dZ$ is symmetric.
    Lastly, we note
    that, for every~$x'$,
    $\PZZ(x, x'')$ contains all the zigzag paths from $x$ to $x''$ passing through $x'$.
    Thus,
    \begin{align*}
        \dZ(x, x'')&= \inf_{\gamma\in\PZZ(x, x'')}\lZ(\gamma)\\
        & \le \inf _{\gamma_1\in\PZZ(x, x'), \gamma_2\in \PZZ(x', x'')} (\lZ(\gamma_1)+\lZ(\gamma_2))\\
        & =\inf _{\gamma_1\in\PZZ(x, x')} \lZ(\gamma_1)+ \inf _{\gamma_2\in\PZZ(x', x'')} \lZ(\gamma_2)\\
        &= \dZ(x, x')+\dZ(x', x'')
    \end{align*}
    showing the triangle inequality for $\dZ$.
\end{proof}

Observe that the zigzag distance is a metric (i.e., every two elements have
finite distance) when the d-space $\X$ is zigzag connected.

\begin{lemma}\label{lem:continuity-directed-wrt-zz}
    Let $(X,d)$ be a metric space together with the directed structure~$\vec P(X)$.
    Every $\gamma\in\PX$
    is continuous with respect to~$\dZ$.
\end{lemma}
\begin{proof}
    First, note that any d-path $\gamma\colon I\to X$ is uniformly continuous with respect to $d$ because $I$ is compact.
    In other words, for any $\varepsilon>0$, there exists~$\delta(\varepsilon)>0$ such that
    \[|s-t| \leq \delta \text{ implies } d(\gamma(s),\gamma(t)) \leq \varepsilon .
    \]
    Let $\varepsilon> 0$ and consider a d-path $\gamma \in \vec P(X)$.
    Let $s,t \in [0,1]$.
    By definition of zigzag distance, for any zigzag path $\xi = (\xi_i)_{i=1}^n
    \in \vec P_{zz}(\gamma(s),\gamma(t))$, we have that
    \begin{equation}\label{eq:d_zz_less_length}
        \dZ(\gamma(s), \gamma(t)) \leq \sum_{i=1}^n \len(\xi_i)
    \end{equation}
    By definition of the length function, $\len$, for any $i$ there exists a partition of $[0,1]$ such that
    \begin{equation}\label{eq:length_sub}
        \len(\xi_i) \leq \sum_{k} d(\gamma(t^{(i)}_k), \gamma (t^{(i)}_{k+1}))+ \frac{\varepsilon}{2n}.
    \end{equation}
    After re-indexing, Equations (\ref{eq:d_zz_less_length}) and (\ref{eq:length_sub}) yield
    \begin{equation} \label{eq:d_zz_Inequality_d}
        \dZ(\gamma(s),\gamma(t)) \leq \sum_{j\in J}d(\gamma(t_j),\gamma(t_{j+1})) + \frac{\varepsilon}{2}.
    \end{equation}
    Inequality (\ref{eq:d_zz_Inequality_d}) still holds under refinements of the partition.
    Hence, without loss of generality, we assume that the partitions of $[0,1]$ were chosen such
    that for any~$j\in J$, we have that $|t_j - t_{j+1}|\leq
    \delta(\frac{\varepsilon}{2|J|})$.
    Then, for any $\varepsilon>0$, we can always find a $\delta>0$ such that
    provided  $|s-t|\leq \delta$, we have that
    \begin{equation*}
        \dZ(\gamma(s), \gamma(t)) \leq \frac{\varepsilon}{2}+\frac{\varepsilon}{2} = \varepsilon.
    \end{equation*}
\end{proof}

\begin{remark}\label{rk:d-metric-space}

    There are now several structures on a metric space $(X,d)$ with a directed
    structure $\PX$, namely the initial metric~$d$, the zigzag metric $\dZ$, and the directed structure $\PX$.
    By Lemma 2, $\PX$ is a directed structure on the metric space $(X,\dZ)$ and
    these are the structures we focus on.
    In particular, we require d-maps to be continuous with respect to the topology induced by $\dZ$.
\end{remark}

\begin{definition}\label{def:directed-ms}
The pair $(\X, \dZ)$ consisting of the (extended) metric space $(X,\dZ)$ together with the directed structure $\PX$ is called an \define{(extended) directed metric space}.
\end{definition}

Given a metric space $(X, d)$, define the intrinsic
metric $d_I \colon X \times X \to \R$ by
\begin{equation}
    d_I (x,x') = \inf_{\gamma \in P(x,x')} \len(\gamma),
\end{equation}
where $P(x,x')$ is the set of all paths starting at $x$ and ending at $x'$ and
$\len$ is defined as in \cref{eqn:length}.  Then
$(X,d_I)$ is the intrinsic metric space associated with the metric $(X,d)$.
If $d_I=d$, then the metric space~$(X, d)$ is called a \textbf{length space}.
We then observe the following:

\begin{proposition}\label{prop:dz-is-d}
    Let $(X, d)$ be a length space, and consider the trivial
    d-structure~$\vec P(X)$ (containing all the continuous paths).
    Then, $\dZ=d$.
\end{proposition}
\begin{proof}
    Consider $x, x'\in X$. Then,
    \[
        \dZ(x, x') = \inf_{\gamma\in \vec P_{zz}(x,x')} \len(\gamma)
        = \inf_{\gamma\in \vec P(x,x')} \len(\gamma)
        = d(x,x').
    \]

\end{proof}

Given a directed metric space $(\vec X, \dZ)$, we can construct a new directed
metric space $(\vec X, (\dZ)_{zz})$ where the d-paths are continuous with
respect to $\dZ$; see \cref{lem:continuity-directed-wrt-zz}.
This
construction can be iterated, and we show in the following proposition that
it is idempotent.

\begin{proposition}\label{prop:idempotency}
    Let $(\X, \dZ)$ be a directed metric space.
    Then, $(\dZ)_{zz}=\dZ$.
\end{proposition}
\begin{proof}
    By \cref{lem:continuity-directed-wrt-zz}, a directed path $\gamma\in \PX$ is
    continuous with respect to~$\dZ$.
    Now, it suffices to show that $\lZ$ is the same with respect to the initial
    metric~$d$ and the associated zigzag distance because the d-structure is the
    same.
    To differentiate the lengths of a zigzag path with respect to the two
    metrics, we denote these lengths by $\lZ$ and $\lZ^{\dZ}$, respectively.
    The \mbox{inequality~$\lZ\le \lZ^{\dZ}$} is a consequence of~$d\le \dZ$; see~\cref{lem_zigzag_bigger_d}.
    On the other hand, consider the zigzag path $\gamma=(\gamma_i)_{i=1}^m$.
    We have that
    \[
        \begin{aligned}
            \lZ^{\dZ}(\gamma)
            &=\sum\limits_{i=1}^{m}\sup_{t^i_j}\sum_j \dZ(\gamma_i(t^i_j),\gamma_i(t^i_{j+1}))\\
            &=\sum\limits_{i=1}^{m}\sup_{t^i_j}\sum_j\inf\limits_{\sigma_i^j\in \vec{P}_{zz}(\gamma_i(t^i_j),\gamma_i(t^i_{j+1}))} \lZ(\sigma_i^j)\\
            &\leq \sum\limits_{i=1}^{m}\sup_{t^i_j}\sum_j \len(\gamma_i^j)
            \\
            & = \sum\limits_{i=1}^{m} \len(\gamma_i)\\
            &=\len_{zz}(\gamma),
        \end{aligned}
    \]
    where $\gamma_i^j$, appearing in the third row, is the subpath of $\gamma_i$
    with $\gamma(t_j)=\gamma_i^j(0)$ and $\gamma(t_{j+1})=\gamma_i^j(1)$.
    Note that the inequality $\lZ^{\dZ}(\gamma)\le \len_{zz}(\gamma)$ implies that
    if $\gamma$ is rectifiable in (X, d) then it is also rectifiable in $(X,
    \dZ)$.
\end{proof}

The assumption that the zigzag distance $\dZ$ is induced from a metric $d$ on a
topological space $X$ is a sufficient condition ensuring that $\dZ$ is an
extended metric.
There are several
other notions of length of a path different from the one in
\cref{def: zigzag_length}.
For instance,
the length may be viewed as
a positive map~$\ell \colon \mathcal{P}(X)\to [0, \infty)$, where
$\mathcal{P}(X)$ is a subset of all paths on $X$ satisfying some axioms
that differ depending on the setting; see,
for instance,~\cite[Def.~1.3]{gromov1999metric}
or~\cite[Sec.~2.1]{mennucci2013asymmetric}.
The following example shows that $\dZ$ need not be an extended metric if we
choose one of these general definitions of length.
In particular, it addresses the importance of starting with a metric space $(X, d)$.

\begin{example}[Open Book]\label{ex:open-book}
    Let $X$ be a topological space and $S\subseteq X$ be a closed subspace of $X$.
    Given $a,b\in S\subseteq X$, with $a\neq b$, and $n\in\mathbb{N}$, define
    $\PX^{(n)}=\langle \gamma_n\rangle=\{c_x\}_{x\in X}\cup \{\gamma_n\circ h\mid h:I\rightarrow I \mbox{ is nondecreasing}\}$,
    where $c_x$ is the constant path at $x\in X$
    and $\gamma_n:
    [0,1]\rightarrow X$ is a simple curve in $X$ such that $\gamma_n(0)=a$
    and~$\gamma_n(1)=b$. Consider the length map~$\ell$
    defined by $\ell(c_x)=0$ and
    $\ell(\gamma_n\circ h)=\frac{1}{n}(t'-t)$ where $[t,t']=h(I)$.
    This implies that $\gamma_n:I\to X$ is parametrized by constant speed $\frac{1}{n}$.

    Then, $(X,\PX^{(n)})$ is a d-space with the induced zigzag metric
    \begin{equation*}
        d^{(n)}(x,x') =
        \begin{cases}
            \ell(\gamma_n\circ h), & \text{if there is $h$  such that} \gamma_n\circ h(0),\gamma_n\circ h(1)\in \{x,x'\},\\
            0,                    & \text{if } x=x',\\
            \infty,               & \text{otherwise. }
        \end{cases}
    \end{equation*}
    Given $n\in\mathbb{N}$, $(S,\P(S))$ is a directed subspace of $(X,\PX^{(n)})$,
    where $\vec{P}(S)=\{c_x\}_{x\in S}$.
    Let $(Y,\PY)=\bigsqcup\limits_{n\in\mathbb{N}}(X, \PX^{(n)})/_\sim$, where
    for every $x\in (X, \PX^{(n)})$ and $x'\in (X, \PX^{(m)})$, for
    some $m\neq n$, $x\sim x'$ iff $x=x'$ and $x,x'\in S$.
    Explicitly, $\PY=\{c_x\}_{x\in Y}\cup \PX^{\infty}$, where $\PX^{\infty}$
    is the smallest set containing all possible concatenations, subpaths, or
    reparameterizations of~\mbox{$\{\gamma_i\mid i\in\mathbb{N}\}$}.
    Equip~$(Y,\PY)$ with the induced pseudo-metric $d_{zz}$.
    Explicitly,
    for all points~$x\in(X,\PX^{(n)})$ and $x'\in(X,\PX^{(m)})$,\\
    \scalebox{0.9}{\parbox{\textwidth}{%
        $$d_{zz}(x, x') = \begin{cases} \min\left\{\begin{array}{l}
            d^{(n)}(x,x'),d^{(n)}(a,x)+d^{(n)}(x',b), \\
            d^{(n)}(a,x')+d^{(n)}(x,b)
        \end{array}\right\} , &
        \begin{tabular}{l}
            if $m=n$,\\
            $x,x'\notin\{a,b\}$;
        \end{tabular}\\
        \min\left\{\begin{array}{l}
            d^{(n)}(a,x)+d^{(m)}(a,x'), d^{(n)}(x,b)+d^{(m)}(a,x'),\\
            d^{(n)}(a,x)+d^{(m)}(x',b), d^{(n)}(x,b)+d^{(m)}(x',b)
        \end{array}\right\} , &
        \begin{tabular}{l}
            if $m\neq n$,\\
            $x,x'\notin\{a,b\}$;
        \end{tabular}\\
        \min\left\{\begin{array}{l}
            d^{(m)}(a,x'), d^{(m)}(x',b)
        \end{array}\right\} , & \begin{tabular}{l}
            if $x\in\{a,b\}$,\\
            $x'\notin\{a,b\}$;
        \end{tabular}\\
        \min\left\{\begin{array}{l}
            d^{(n)}(a,x), d^{(n)}(x,b)
        \end{array}\right\} , & \begin{tabular}{l}
            if $x'\in\{a,b\}$,\\
            $x\notin\{a,b\}$;
        \end{tabular}\\
        0, & \mbox{if } x,x'\in \{a,b\}.\\
\end{cases}$$}}

Note that $d_{zz}(a,b)=0$ because the infimum of the length of the zigzag
paths connecting $a$ and $b$ is $\inf\limits_{n\in\mathbb{N}}\frac{1}{n}=0$,
despite the assumption of $a\neq b$. This suggests that even if the zigzag
distance could be defined also in the absence of a metric in the underlying
space, it would not be an extended metric in general, but an extended pseudo-metric.
\end{example}

\begin{example}[Non-Equivalent Topology]\label{ex:non-eq-top}
    By \cref{lem:continuity-directed-wrt-zz} the d-paths are continuous
    with respect to the topology induced by the zigzag distance but, in
    general, this topology is not equivalent to the topology of the underlying metric space.

    To illustrate this, take~$(X,d)$ to be the Euclidean plane. Let the set of
    d-paths~$\PX$ in the d-space
    $(X,\PX)$ be generated by
    \begin{align*}
        \{ \gamma_x \colon I \to X,\ \gamma_x(t) = t \cdot x \mid x \in X \}.
    \end{align*}
    The zigzag metric,~$\dZ$, that these choices induce, is known by many names,
    including post office, French metro, British rail and SNCF\footnote{SNCF is
    the French railway, the \emph{Soci\'et\'e nationale des chemins de fer
    français}
    (National Company of French Railways).} metric.
    Fix a radius~$R>0$ and take a point~$x \in X$ for which~$d(0,x) > R$. Notice that, for
    any~\mbox{$r\geq R$},~$B_{zz}(x,R) \subset B(x,r)$.
    However,  $B_{zz} (x, R)$ is a segment of diameter~$2R$ along the line connecting $0$ to $x$, so
    there exists no~$r' > 0$ such that
    \begin{align*}
        B(x,r') \subseteq B_{zz}(x,R),
    \end{align*}
    where~$B$ is a ball in the Euclidean metric~$d$, and~$B_{zz}$ a ball in the induced
    zigzag distance~$\dZ$.
    Not all the paths that are continuous with respect to $d$ are still
    continuous with respect to $\dZ$.
    For example, the  path $\gamma(t)=(\cos(t), \sin(t))$ is not continuous with
    respect to the zigzag distance.
\end{example}

In contrast to Example~\ref{ex:non-eq-top}, we show an interesting example where
the topologies induced by a metric and the induced zigzag distance are
equivalent.

\begin{example}[Directed Flat Torus]
    Consider the rectifiable d-space $(\vec{I}^2, \dZ)$, where the
    d-structure is induced by the product order on the unit square~$I^2$ and $\dZ$ is induced
    by the~$L^2$-distance, and the flat torus $\mathbb{T}=I^2/\sim$, with
    $(0,y)\sim (1,y)$ and $(x,0)\sim (x,1)$.
    As described in \cref{subsec:operations}, $I^2/\sim$ inherits a rectifiable
    d-structure which is the smallest one making the quotient map $\pi\colon
    I^2\to I^2/\sim$ a d-map.
    In general, when passing to the quotient, we obtain a pseudo-metric.
    However, for $\mathbb{T}$, the distance induced by the $L^2$-distance on
    $I^2$ is, indeed, a metric.
    The resulting d-space, $\vec{\mathbb{T}}$, is called the directed flat
    torus.
    Explicitly, d-paths in~$\mathbb{T}$ are concatenations of~$\pi\circ \gamma$,
    with $\gamma \in \P(I^2)$, and reparametrizations thereof.

    In the zigzag metric, a ball on the flat torus with a small radius
    (\cref{subfig:torus_small_ball}) looks like a combination of the balls
    in~$L^1$ and~$L^2$ norms. When we increase the radius sufficiently, the ball
    starts
    overlapping with itself, as shown in
    \cref{subfig:torus_big_ball}.
    \begin{figure}[ht]
        \begin{subfigure}{.49\textwidth}
            \begin{tikzpicture}[scale=0.55,
	dot/.style = {circle, fill=black, inner sep=0pt, minimum size=5pt, node contents={}}
	]
	
	\coordinate (a) at (-5,-5);
	\coordinate (b) at (5,-5);
	\coordinate (c) at (5,5);
	\coordinate (d) at (-5,5);
	\coordinate (center) at (0, 0);
	
	\draw (a) -- (b) -- (c) -- (d) -- (a);
	\path (0, -5) node {$\gg$};
	\path (0, 5) node {$\gg$};
	\node[label={[rotate=90]center:$>$}] at (-5,0) {};
	\node[label={[rotate=90]center:$>$}] at (5,0) {};

	\fill[green!20!white] (0,0) -- (2cm,0mm) arc (0:90:2cm) -- (0,0);
	\fill[green!20!white] (0,0) -- (-2cm,0mm) arc (180:270:2cm) -- (0,0);
	\fill[green!20!white] (-2,0) -- (0,-2) -- (2,0) -- (0,2) -- cycle;

	\draw[step=0.5cm, opacity=0.3] (-5,-5) grid (5, 5);
	\path (center) node[dot] {};
	
\end{tikzpicture}

            \caption{Ball with a small radius.}
            \label{subfig:torus_small_ball}
        \end{subfigure}
        \hfill
        \begin{subfigure}{.49\textwidth}
\begin{tikzpicture}[scale=0.55,
	dot/.style = {circle, fill=black, inner sep=0pt, minimum size=5pt, node contents={}}
	]

	\coordinate (a) at (-5,-5);
	\coordinate (b) at (5,-5);
	\coordinate (c) at (5,5);
	\coordinate (d) at (-5,5);
	\coordinate (center) at (0, 0);
	
	\draw (a) -- (b) -- (c) -- (d) -- (a);
	\path (0, -5) node {$\gg$};
	\path (0, 5) node {$\gg$};
	\node[label={[rotate=90]center:$>$}] at (-5,0) {};
	\node[label={[rotate=90]center:$>$}] at (5,0) {};
	
	\clip (-5,-5) rectangle (5,5);
	\fill[green!20!white, opacity=0.8] (0,0) -- (6cm,0mm) arc (0:90:6cm) -- (0,0);
	\fill[green!20!white, opacity=0.8] (0,0) -- (-6cm,0mm) arc (180:270:6cm) -- (0,0);
	\fill[green!20!white, opacity=0.8] (0,-6) -- (6,0) -- (0,0) -- cycle;
	\fill[green!20!white, opacity=0.8] (-6,0) -- (0,6) -- (0,0) -- cycle;
	
	\fill[blue!20!white, opacity=0.8] (0,-10) -- (6,-10) arc (0:90:6cm) -- (0,-10);
	\fill[blue!20!white, opacity=0.8] (-1,-5) -- (0,-5) -- (0,-4) -- cycle;
	\fill[blue!20!white, opacity=0.8] (-10,0) -- (-4,0) arc (0:90:6cm) -- (-10,0);
	\fill[blue!20!white, opacity=0.8] (-5,0) -- (-5,-1) -- (-4,0) -- cycle;
	\fill[blue!20!white, opacity=0.8] (0,10) -- (-6,10) arc (180:270:6cm) -- (0,10);
	\fill[blue!20!white, opacity=0.8] (1,5) -- (0,5) -- (0,4) -- cycle;
	\fill[blue!20!white, opacity=0.8] (10,0) -- (4,0) arc (180:270:6cm) -- (10,0);
	\fill[blue!20!white, opacity=0.8] (4,0) -- (5,0) -- (5,1) -- cycle;
	
	\draw[step=0.5cm, opacity=0.3] (-5,-5) grid (5,5);
	\path (center) node[dot] {};
	
	\path (0, -5) node {$\gg$};
	\path (0, 5) node {$\gg$};
	\node[label={[rotate=90]center:$>$}] at (-5,0) {};
	\node[label={[rotate=90]center:$>$}] at (5,0) {};
	
\end{tikzpicture}

            \caption{Ball with a big radius.}
            \label{subfig:torus_big_ball}
        \end{subfigure}
        \caption{Balls in the zigzag metric on a directed flat torus.}
        \label{fig:balls_on_torus}
    \end{figure}
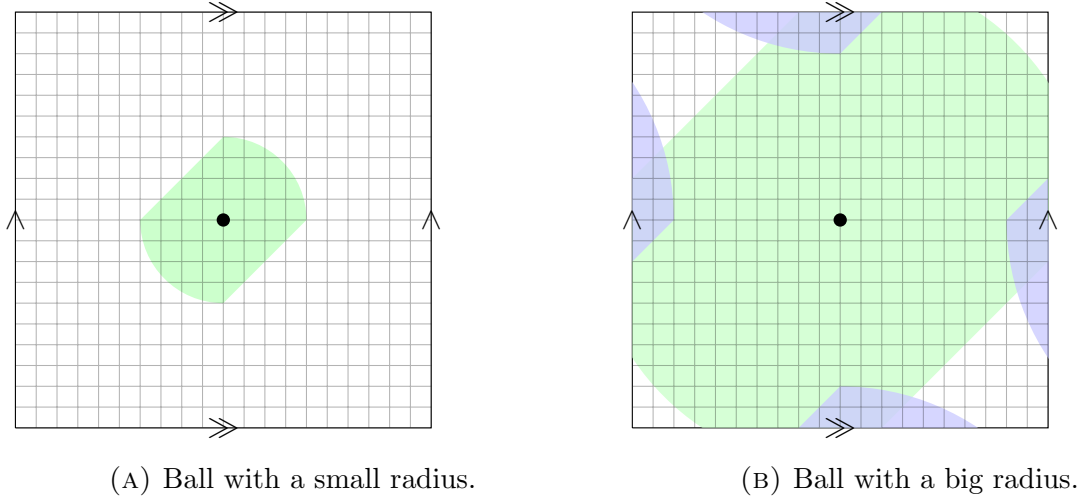
    In this example the topology given by the zigzag metric is
    equivalent to the topology given
    by the metric inherited from the $L^2$-distance via the quotient.
    This is because
    \begin{align*}
        B\Big(\frac{r \sqrt{2}}{2}\Big) \subset B_{zz}(r) \subset B(r)
    \end{align*}
    for any~$r>0$, where~$B$ denotes balls in the $L^2$-distance,
    and~$B_{zz}$ denotes balls in the zigzag metric induced by it.
\end{example}

\begin{example}[Directed Weighted Graphs]
    Given a directed graph as in \cref{ex:dirgraph} together with a
    positive weight function~$W \colon E \to \R_{>0}$, we further
    modify the construction of~$X_G$ to incorporate the weight of an
    edge as its length. The gluing map for each
    edge~$e \in E$ is defined as~$\varphi_e \colon
    \partial
    \overrightarrow{[0, W(e)]} \to X^0_G$ with~$\varphi_e(0)=s(e)$ and~$\varphi_e
    (W(e)) = t(e)$, where~$s$ and~$t$ map each edge to its source and
    target respectively. This enables one to define a metric
    structure on~$X_G$, much like in the context of metric
    graphs~\cite{aanjaneya2011metric,mugnolo2019actually}.
    Each path~$\gamma \colon I \to X_G$ is split
    into subpaths~$\gamma_i$ so that, for every $i$, $\gamma_i([0,1])$
    is contained entirely within one {one-cell}, $e_i$. Then, the length of
    such a path is the sum~$\sum_i \len_i(\gamma_i)$, where~$\len_i$ is the
    length within the corresponding cell,~$e_i$. Defining a distance
    between two points as the infimum of the paths connecting them
    gives a metric, proving this is an easy exercise we leave to the reader.
\end{example}

\subsection{Directed Gromov--Hausdorff Distance}\label{sub:directed-GH-distance}
From now on, unless stated otherwise, all directed spaces
$\X = (X, \PX)$ are endowed with the zigzag metric, $\dZ^X$, induced from the
underlying metric space~$(X,d^{X})$; recall \cref{rk:d-metric-space}.

In order to define the Gromov--Hausdorff distance for d-spaces, we first need a
notion of a directed isometry.

\begin{definition}\label{def:directed-isometry}
    Consider two directed metric spaces $(\X,\dZ^{X})$ and $(\Y,\dZ^{Y})$.
    A d-map $\dF \colon \X\to\Y$ is called a \define{d-isometry} if for any $x,
    x'$ in $X$,
    \[
        \dZ^{X}(x, x')=\dZ^{Y}(\dF(x), \dF(x')).
    \]
    Two d-spaces $\X$ and $\Y$ are \define{d-isometric} if and only if there is a bijective
    d-isometry $\dF\colon \X\to \Y$ such that its inverse (as a function
    $F^{-1}\colon Y\to X$) is a~d-map.
\end{definition}

Let $\X$ and $\Y$ be d-subspaces of the directed metric space $(\Z, \dZ)$.
The
\define{directed Hausdorff distance} of $\X$ and $\Y$ is defined by
\begin{align}\label{eqn:H}
    \dhaus(\X, \Y)
    =\haus((X, \dZ), (Y, \dZ)).
\end{align}

\begin{definition}\label{def:d_GH1}
 We define the \define{directed Gromov--Hausdorff distance} between two
    directed metric spaces
    $(\X, \dZ^X)$ and $(\Y, \dZ^Y)$ by:
    \begin{align}\label{eqn:directed-GH}
        \dGH(\X,\Y) = \inf_{\dF, \dG} \dhaus(\dF(\X), \dG(\Y)),
    \end{align}
    where $\dF\colon \X \to \Z$ and $\dG\colon \Y \to
    \Z$ are d-isometries into some directed metric space $(\Z, \dZ^Z)$.
\end{definition}

Imposing that the isometries between directed spaces in the definition of the
Gromov--Hausdorff distance are d-maps is not restrictive, as shown by the following result.
The Gromov--Hausdorff distance, in fact, depends only on the zigzag metric
structure induced by the d-structures.

\begin{theorem}\label{thm:gh-vs-dgh}
    Let~$(\X, \dZ^X)$ and~$(\Y, \dZ^Y)$ be directed metric spaces.
    Then,
    \begin{equation*}
        \dGH((\X, \dZ^X), (\Y, \dZ^Y)) = \GH((X,\dZ^{X}), (Y, \dZ^{Y})).
    \end{equation*}
\end{theorem}
\begin{proof}
    We observe that the set of d-isometries from the d-spaces $(\X, \dZ^X)$ and~$(\Y,
    \dZ^Y)$ to a d-space~$(\Z, \dZ^Z)$ is contained in the set of isometries from
    $(X, \dZ^X)$ and~$(Y, \dZ^Y)$ to $(Z, d^Z)$. Thus,~$\GH((X, \dZ^X), (Y, \dZ^Y))\le \dGH((\X, \dZ^{X}), (\Y,\dZ^{Y})) $.
    Let us therefore prove the opposite inequality.
    By definition of the classical Gromov--Hausdorff distance, for any
    $\varepsilon>0$, there exists a metric space $(Z,d^{Z})$ and two
    isometries~$F: (X, \dZ^{X}) \hookrightarrow (Z,d^{Z})$ and~$G: (Y, \dZ^{Y}) \hookrightarrow (Z, d^{Z})$ such that
    \begin{equation*}
        \haus(F(X), G(Y)) \leq \GH((X,\dZ^{X}), (Y, \dZ^{Y})) + \varepsilon.
    \end{equation*}
    By \emph{Kuratowski's embedding theorem}, every metric space isometrically
    embeds into a Banach space.
    Then, without loss of generality, we can always choose $Z$ to be a Banach space.
    Let
    $\P(Z) = \langle \lbrace F \circ \gamma : \gamma  \in \P(X) \rbrace
    \cup \lbrace G \circ \gamma : \gamma \in \P(Y) \rbrace \cup \vec L\rangle$,
    where $\vec L$ is the d-structure generated
    by~$\lbrace (1-t)F(x)+tG(y) : x \in X, y\in Y \rbrace$ of line segments connecting points in $F(X)$ and $G(Y)$ in $Z$.
    Then,~$(Z, \P(Z))$ is a d-space. Note that for any $x_1,x_2 \in X$,
    \begin{equation*}
        \dZ^{Z}(F(x_1),F(x_2)) = \inf_{\gamma \in \PZZ^{Z}(F(x_1), F(x_2))} \text{len}(\gamma) \leq  \inf_{\gamma\in \PZZ^{X}(x_1,x_2)} \text{len}(\gamma) = \dZ^{X}(x_1,x_2),
    \end{equation*}
  where the inequality follows from the fact that $F$ is an isometry.
  Moreover, because
  \begin{equation*}
      \dZ^{X}(x_1,x_2)=d^{Z}(F(x_1),F(x_2)) \leq \dZ^{Z}(F(x_1),F(x_2)),
  \end{equation*}
  we have that
  $F\colon (X, \dZ^X)\to (Z, \dZ^Z)$ is an isometry.
  Similarly, one can show that
  $G\colon (Y, \dZ^Y)\to (Z, \dZ^Z)$ is an isometry.
    Because $F$ and $G$ map d-paths in~$\PX$ and~$\PY$, respectively, to d-paths in
    $\P(Z)$, we conclude that $F\colon \X \rightarrow \Z$ and~$G\colon \Y \rightarrow \Z$ are d-isometries.

We now show that for any $x \in X, y\in Y $ and for any $\gamma \in \PZZ(x,y)$, $\dZ^{Z}(x,y) = d^{Z}(x,y)$.
On the one hand,
\[
\dZ^Z(F(x), G(y))\le \len((1-t)F(x)+tG(y))=d^Z(F(x), G(y)),
\]
where the equality holds because the length of the line segments in a Banach space equals the distance between their end-points.
On the other hand, $d^Z\le \dZ^Z$, by \cref{prop_zigzag_metric}.
This implies that
\begin{equation*}
    \haus((F(X),\dZ^{Z}),(G(Y), \dZ^{Z}) ) = \haus((F(X),d^{Z}),(G(Y), d^{Z})).
\end{equation*}
Finally,
\begin{equation*}
    \begin{split}
         \dGH((\X, \dZ^X), (\Y, \dZ^Y)) &\leq \haus((F(X),\dZ^{Z}),(G(Y), \dZ^{Z}) ) \\
         & \leq \GH((X,\dZ^{X}), (Y, \dZ^{Y})) + \varepsilon.
    \end{split}
\end{equation*}
The result follows by sending $\varepsilon$ to zero.
\end{proof}

\begin{remark}\label{dirGH-length}
    By \cref{thm:gh-vs-dgh} and \cref{prop:dz-is-d}, if $(X, d^X)$ and $(Y,
    d^Y)$ are length spaces endowed with the trivial d-structure, then the
    directed Gromov-Hausdorff distance between them coincides with the classical
    Gromov-Hausdorff distance.
\end{remark}

Recall that the Gromov--Hausdorff distance is a metric on the space of isometry
classes of compact metric
spaces. 
If we consider the isometry classes with respect to zigzag metrics, we deduce from \thmref{gh-vs-dgh} the following result.

\begin{corollary}\label{cor: d-GH metric}
    The directed Gromov--Hausdorff distance is a metric on the space of isometry
    classes of compact directed metric spaces.
\end{corollary}

Note that \cref{cor: d-GH metric} requires compactness of the directed metric
space~$(\X, \dZ
^X)$.
This is a stronger requirement than that of compactness of $(X, d)$.
In fact, compactness of $(\X, \dZ^X)$ implies compactness of $(X, d)$, as $d$
induces a coarser topology than that induced by $\dZ$, but not vice versa.
As an example of a compact space $(X, d)$ with non-compact associated directed metric space $(\X, \dZ)$, consider $(I^2, L^2)$ with the d-structure generated by all continuous paths contained in $[0,1)\times [0,1]\cup\{(1,1)\}$ and the straight-line path from $(1,1)$ to $(1,0)$.
The Cauchy sequence $a_n=(1-\frac{1}{n}, 0)$ does not converge.
Thus, the directed metric space $(\vec I^2, L^2_{zz})$ is not complete.

Further consequences of \thmref{gh-vs-dgh} for directed metric spaces are
summarized in the following two corollaries.

First, let us observe that if two metric spaces endowed with a d-structure
induce the same zigzag distance, then the directed Gromov--Hausdorff distance
between them is zero.
This might happen in the cases expressed in \cref{rem:samedzz}, for example.
In particular, in the following corollary, we use observation from
\cref{ex:reverse} and \cref{rem:samedzz} that $\vec{X}$ and $\Xrev$ induce the
same zigzag distance~$\dZ^X$.
\begin{corollary} \label{prop:GH_X_rX}
    For any directed metric space $(\X, \dZ^X)$, we have $\dGH(\X, \Xrev)=0$.
\end{corollary}

Combining \thmref{gh-vs-dgh} and Proposition \ref{prop:GH_properties}, we obtain
the following results.

\begin{corollary}\label{cor:prop-dirGH}
    Let $(\X, \dZ^X)$ and $(\Y, \dZ^Y)$ be directed metric spaces.Then,
    \begin{enumerate}
        \item $\dGH(\X, \Y) \le \frac{1}{2} \max \lbrace \diam(\X, \dZ^X) , \diam(\Y,\dZ^Y)\rbrace$.
        \item $\dGH(\X, \Y) = \frac{1}{2} \diam(\Y, \dZ^Y)$, if $X=\{x_0\}$.
        \item\label{item:dirGH-lowerbd} $\dGH(\X, \Y)\geq \frac{1}{2} |\diam(\X, \dZ^X)-\diam(\Y, \dZ^Y)|$, if $\diam(\X, \dZ^X)<\infty$.
    \end{enumerate}
\end{corollary}

We note that, despite the result stated in \cref{cor: d-GH metric}, the directed
Gromov--Hausdorff distance is not a metric on the space of compact directed
metric spaces up to d-isometry. Indeed, we show next an example where $\dGH$ is degenerate.

\begin{example}\label{ex:dGH_not_metric}
    Let~$X = [-1,1]$ and define the maps~$\gamma_1, \gamma_2\colon I \to X$
    as~\mbox{$\gamma_1(t)= t$} and $\gamma_2(t)= -t$.
    Endow $X$ with the path space~$\PX = \langle \gamma_1, \gamma_2 \rangle$ and
    the Euclidean metric.
    Denote the resulting d-space as~$\X$ and its reverse space as~$\Xrev$;
    see~\cref{ex:reverse}.
    \begin{center}
\begin{tikzpicture}[scale=1.6,
	dot/.style = {circle, fill=black, inner sep=0pt, minimum size=5pt, node contents={}},
	midarrow/.style = {decoration={markings, mark=at position 0.5 with {\arrow{>}}},
		postaction={decorate}}
	]
	
	\coordinate (a) at (-3,0);
	\coordinate (b) at (-2,0);
	\coordinate (c) at (-1,0);
	\coordinate (c') at (3,0);
	\coordinate (b') at (2,0);
	\coordinate (a') at (1,0);
	
	\path (a) node[dot,label=above:$-1$];
	\path (b) node[dot,label=above:$0$];
	\path (c) node[dot,label=above:$1$];
	\path (a') node[dot,label=above:$-1$];
	\path (b') node[dot,label=above:$0$];
	\path (c') node[dot,label=above:$1$];
	
	\draw[midarrow] (b) -- (a) node[midway, above=1.5ex]{$\gamma_2$};
	\draw[midarrow] (b) -- (c) node[midway, above=1.5ex]{$\gamma_1$};
	\draw[midarrow] (a') -- (b') node[midway, above=1.5ex]{$\gamma_2^\ast$};
	\draw[midarrow] (c') -- (b') node[midway, above=1.5ex]{$\gamma_1^\ast$};
	\path (a) node[left=1.5em]{$\X:$};
	\path (a') node[left=1.5em]{$\X^\ast:$};
\end{tikzpicture}
    \end{center}
    By \cref{prop:GH_X_rX},
    $\dGH \big ( (\X, \dZ^X),(\Xrev, \dZ^X) \big ) = 0$. \par

    For any d-map $\dF \colon \X \to \Xrev$, the compositions $F
    \circ \gamma_1$ and $F \circ \gamma_2$ must be d-paths in
    $\PX^\ast$.
    Further, because $F \circ \gamma_1(0) = F(0) = F \circ
    \gamma_2(0)$, one of the following holds:
    \begin{itemize}
        \item $F \circ \gamma_1 (I), F \circ \gamma_2 (I) \subseteq \gamma_1^\ast(I) = [0, 1]$,
        \item $F\circ \gamma_1 (I), F \circ \gamma_2 (I)\subseteq \gamma_2^\ast (I) = [-1, 0]$.
    \end{itemize}
    As a consequence, there is no d-isometry between $\X$ and $\Xrev$.
    Note, however, that both $\X$ and $\X^*$ with the zigzag distance are
    isometric to $X$ with the Euclidean~metric.
\end{example}

Now, we show an example where the distance between $(X, d^X)$
and $(X, \dZ^X)$ is non-zero.

\begin{example}\label{ex:GH_Xd_Xdzz} Let $X = I^2$ and let $d^X$ be the
    Euclidean metric. Let $(X, \PX)$ be the d-space induced by the product order
    on $I^2$. Then, we show that
    \[
        \GH\left((X, d^X), (X, \dZ^{X}) \right) = 1 - \frac{\sqrt{2}}{2}.
    \]
    Indeed, from \cref{prop:GH_properties}.4 we know \[\GH\left((X, d^X), (X, \dZ^{X}) \right)  \geq \frac{1}{2}\left|\diam(X, d^X) -  \diam(X, \dZ^X)\right| = \frac{1}{2}\left(2 - \sqrt{2}\right).\]

    Using \cref{eqn:kalton}, we obtain this lower bound by letting~$(X, \dZ^X)
    \xrightarrow{\varphi} (X, d^X)$ and $(X, \dZ^X) \xleftarrow{\psi} (X, d^X)$
    be the identity maps. In this case, only noncomparable points contribute
    to the (co)distortion, which is given by maximizing the following:~$|L^1(x, x') -
    L^2(i(x), i(x'))|$.
    This maximum is obtained by antidiagonal points. Thus
    \[
        \dis(\varphi) = \dis(\psi) = \codis(\varphi, \psi) = 2 - \sqrt{2},
    \]
    proving that $\GH\left((X, d^X), (X, \dZ^{X}) \right) = 1 - \frac{\sqrt{2}}{2}$.
\end{example}

\subsection{Distortion and Codistortion of Directed Maps}
\label{sub:distortion-distance}
In this subsection all the directed metric spaces $(\X, \dZ^X)$ are assumed to
be bounded.
Moreover, whenever we assume boundedness for $\X$, it implicitly refers to the
the zigzag distance $\dZ^X$.

The identity from \cref{eqn:kalton} connects the Gromov--Hausdorff
distance with distortion
and codistortion, two quantities associated with maps between spaces.
It is very useful, as it provides bounds
for otherwise notoriously difficult-to-compute Gromov--Hausdorff distance; see,
for example,~\cite{adams2022gromov,lim2021gromov}.
We explore whether a similar connection exists for d-spaces.

\begin{remark}\label{rem:alternate-dgh}
    Observe that, as an immediate consequence of \thmref{gh-vs-dgh}, we obtain
    equivalent definitions of the directed Gromov-Hausdorff distance between directed metric spaces.
    First, as a consequence of \cref{eqn:kalton}, we obtain the following when $\X$ and $\Y$ are bounded:
    \begin{equation}\label{eqn:alternate-dgh-fg}
        \dGH(\X, \Y)
        = \frac{1}{2} \inf_{f, g} \max\{\dis(f), \dis(g), \codis(f, g)\},
    \end{equation}
    where we infimize over (not necessarily continuous) maps $f: (X, \dZ^X) \to
    (Y, \dZ^Y)$ and $g: (Y, \dZ^Y) \to (X, \dZ^X)$.
    Alternatively, as a consequence of \cref{eqn:kalton_bis}, we get a second
    equivalent definition:
    \begin{equation}\label{eqn:alternate-dgh-correspondence}
        \dGH(\X, \Y)
        = \frac{1}{2} \inf_{\mathcal{R}} \sup_{(x, y), (x’, y’) \in \mathcal{R}} |\dZ^X(x, x’) - \dZ^Y(y, y’)|,
    \end{equation}
    where $\mathcal{R}$ ranges over all the correspondences between $X$ and $Y$.
\end{remark}

However, note that Equations~\eqref{eqn:alternate-dgh-fg}
and~\eqref{eqn:alternate-dgh-correspondence} do not utilize the directed
structure of $\X$ and $\Y$. Much like our definition of the directed
Gromov--Hausdorff distance, they only consider the topology induced by the
zigzag distance. In this section, we propose other distances between directed
metric spaces that more naturally account for the path structures.

We note that distortion and codistortion are also defined for d-maps. That
is, if $\dF \colon (\X, \dZ^X) \to (\Y, \dZ^Y)$ is a d-map, we define its
distortion by
\begin{equation}
    \dis(\dF) = \sup_{x, x' \in \X} |\dZ^X(x,x') - \dZ^Y(\dF(x), \dF(x'))|.
\end{equation}
Codistorsion is defined analogously. 

\begin{definition}\label{d-dis-distance}
    The \define{distortion distance}
    between bounded directed metric spaces~$(\X, \dZ^X)$ and $(\Y, \dZ^Y)$ is
    defined as
    \begin{align}
        \ddis(\X, \Y)
        = \frac{1}{2}\inf_{\dF, \dG} \max\{\dis(\dF), \dis(\dG), \codis(\dF, \dG)\},
    \end{align}
    where $\dF\colon \X \to \Y$ and $\dG \colon \Y \to \X$ are d-maps.
\end{definition}

\begin{proposition}\label{lem:distortion_distance_triangle}
Let $\X$, $\Y$, and $\Z$ be bounded directed metric spaces. Then,
\[
\ddis(\X, \Z) \leq \ddis(\X, \Y) + \ddis(\Y, \Z).
\]
\end{proposition}
\begin{proof}
    Consider the compositions
    $\X \xrightarrow{\dF_1} \Y \xrightarrow{\dG_1} \Z \;$
    and~$\X \xleftarrow{\dF_2} \Y \xleftarrow{\dG_2} \Z$.
    Utilizing the triangle inequality yields
    \begin{align*}
   & \sup_{x, x' \in X} |\dZ^X(x,x') - \dZ^Z(\dG_1 \circ \dF_1(x), \dG_1 \circ \dF_1(x') )| \\
   \leq &\sup_{x, x' \in X} |\dZ^X(x,x') - \dZ^Y( \dF_1(x), \dF_1(x') )| + \sup_{y, y' \in Y} | \dZ^Y(y, y' )- \dZ^Z(\dG_1(y), \dG_1(y') ) |
    \end{align*}
    That is,
    \begin{equation}\label{eqn:distortion_traingle}
        \dis(\dG_1 \circ \dF_1) \leq \dis(\dG_1) + \dis(\dF_1)
    \end{equation}
    and a similar argument shows that \cref{eqn:distortion_traingle}
    also holds for $\dis(\dF_2 \circ \dG_2)$
    and $\codis(\dG_1 \circ \dF_1, \dF_2 \circ \dG_2)$.
    Now, consider instead the (larger) set of maps $\vec{\varphi}: \X \to \Z$
    and $\vec{\psi}: \Z \to \X$.~Then,
    \begin{align*}
        \ddis(\X, \Z)
        &= \frac{1}{2}\inf_{\vec{\varphi}, \vec{\psi}} \max \{\dis(\vec{\varphi}), \dis(\vec{\psi}), \codis(\vec{\varphi}, \vec{\psi}) \} \\
        & \leq \frac{1}{2}\inf_{\substack{\dG_1 \circ \dF_1,\\ \dF_2 \circ \dG_2}} \max \{\dis(\dG_1 \circ \dF_1),
        \dis(\dF_2 \circ \dG_2),
        \codis(\dG_1 \circ \dF_1, \dF_2 \circ \dG_2) \} \\
        &\leq \frac{1}{2}\inf_{\dF_1, \dF_2} \max \{\dis(\dF_1), \dis(\dF_2), \codis(\dF_1, \dF_2) \} \\
        &\qquad + \frac{1}{2}\inf_{\dG_1, \dG_2} \max \{\dis(\dG_1), \dis(\dG_2), \codis(\dG_1, \dG_2) \}\\
        & = \ddis(\X, \Y) + \ddis(\Y, \Z).
    \end{align*}
\end{proof}

\begin{theorem}\label{thm:distortion0_iff_disometric}
    Directed maps between bounded directed spaces~$\dF: \X \to \Y$ and~$\dG: \Y \to \X$ for
    which $\dis(\dF) = \dis(\dG) = \codis(\dF,\dG) = 0$ exist if and
    only if the spaces~$\X$
    and~$\Y$ are d-isometric.
\end{theorem}
\begin{proof}
    Suppose such maps exist and prove the forward implication.
    Let~$x, x' \in X$ be such that~$\dF(x) = \dF(x')$.
    Then,~$\dZ^X(x, x') -\dZ^Y(\dF(x), \dF(x')) = \dZ^X(x, x')$,
    and because~$\dis(\dF)=0$, we have
    $\dZ^X(x, x')=0$. Because~$\dZ^X$ is a metric,~$x=x'$.

    A similar argument shows~$\dG$ is injective.

    They are also surjective, which we show for~$\dF$. Suppose there is~$y \in
    Y$ such that~$y \notin \dF (X)$.
    This means that for any~$x \in X$,~$\dZ^X(x, \dG(y)) = \dZ^Y (\dF(x),
    y)>0$, where the
    equality follows from~$\codis(\dF, \dG)=0$. Choose~$x
    = \dG(y) \in X$, this gives~$\dZ^X
    (x, \dG(y)) = \dZ^Y (\dF(x), y)=0$. By contradiction,~$Y = \dF(X)$.

    Moreover, because the codistorsion is $0$,  for any~$y \in Y$,
    \begin{align}
        \dZ^Y (\dF(\dG(y)), y)=\dZ^X(\dG(y), \dG(y))=0,
    \end{align}
    we have that~$\dF \circ \dG = \id_Y$.
    Similarly,~$\dG \circ \dF = \id_X$, which shows $\dF$ and $\dG$ are inverses
    of each other.

    The backwards implication follows easily.
    Let~$\dF \colon \X \to \Y$ be a
    bijective d-isometry and~$\dG \colon \Y \to \X$ its inverse.
    Because both maps
    are d-isometries,~$\dis(\dF) = \dis(\dG)=0$, and because~$\dF$ is a d-isometry,
    \begin{align*}
        d^X(x, \dG(y)) = d^Y(\dF(x), \dF \circ \dG(y)) = d^Y(\dF(x), y)
    \end{align*}
    for all~$x \in X$ and~$y \in Y$. Consequently,~$\codis(\dF, \dG) = 0$.
\end{proof}

\begin{remark}
    Combining \cref{lem:distortion_distance_triangle} and
    \cref{thm:distortion0_iff_disometric} shows that $\ddis$ is an extended
    pseudo-metric
    on the space of d-isometry classes of compact directed metric spaces (symmetry is
    obvious).
    We have not yet shown it is a metric, as we have not shown that $\ddis(\X,
    \Y) = 0$ implies $\X$ and $\Y$ are d-isometric (due to the infimum in the
    definition of $\ddis$).
\end{remark}

\begin{proposition}\label{l:dK_greaterorequal_dGH}
For any bounded 
directed metric spaces $\X$ and~$\Y$,
\[
   \dGH(\X, \Y)\le \ddis(\X, \Y).
\]
\end{proposition}

\begin{proof}
    Recall identity~\eqref{eqn:kalton} and
    \cref{thm:gh-vs-dgh}.
    Because the family of
    pairs~$(\dF\colon \X \to \Y, \dG \colon \Y \to \X)$ of d-maps is
    included in the set of all maps between the two spaces (not necessarily
    directed), the proposition follows.
\end{proof}

The following provides an example in which $\ddis$ is strictly
bigger than~$\dGH$.

\begin{example}\label{ex:dK_greaterorequal_dGH}
    Consider $\X$ and $\Xrev$ as in \cref{ex:dGH_not_metric}.
    We observed there that the image of every d-map $\dF\colon \X\to
    \Xrev$ is either included in $[0,1]$ or in $[-1,0]$.
    As a consequence, $\dZ^X(-1, 1) - \dZ^X(\dF(-1), \dF(1))
    \geq 2 - 1$, and~$\dis(\dF) \geq 1$.
    By symmetry we also have that~$\dis(\dG) \geq 1$, for any d-map~${\dG
    \colon \Xrev \to \X}$. This means~that
    \begin{align} \label{eqn:max_in_sourcesink}
        \max \{\dis(\dF), \dis(\dG), \codis(\dF, \dG)\} \geq 1
    \end{align}
    for any pair of d-maps~$\dF \colon \X \to \Xrev$ and ${\dG \colon
    \Xrev \to \X}$.
    This already implies that~$\ddis(\X, \Xrev)>\GH(\X,
    \Xrev)= 0$; see \cref{ex:dGH_not_metric}.
    However, for this example, we show that
    $\ddis(\X, \Xrev)=\frac{1}{2}$.
    Let us define the d-maps $\dF \colon \X \to \Xrev$ and~$\dG \colon
    \Xrev \to \X$ as follows:
    \begin{align*}
        \dF(t) =
        \begin{cases*}
            -1, & for $t \in [-1, 0]$, \\
            t-1, & for $t \in [0,1]$,
        \end{cases*}
        \qquad
        \dG(t) =
        \begin{cases*}
            t+1, & for $t \in [-1, 0]$, \\
            1, & for $t \in [0,1]$.
        \end{cases*}
    \end{align*}
    To compute the distortion of~$\dF$, first notice that~$\dZ^X(x,
    y) - \dZ^X(\dF(x), \dF(y))$ is at most~$1$ when~${x,
    y \in [-1,0]}$ and it is~$0$ when~$x, y \in [0,1]$.
    Now, assume~$x \in [-1,0]$ and~$y \in [0,1]$. Then,
    \begin{align*}
        \dZ^X(x, y) - \dZ^X(\dF(x), \dF(y))
        = y - x - |-1 - (y -1)| = -x,
    \end{align*}
    which attains a maximum value of~$1$ when~$x = -1$. Thus~$\dis(\dF) = 1$.
    We follow similar steps to show~$\text{dis}(\dG) = 1$ as well.
    Next, let us compute the codistortion of~$\dF$ and~$\dG$.
    By following the definitions of~$\dF$ and~$\dG$, we see that
    \begin{align*}
        |\dZ^{X}(x, \dG(y)) - \dZ^X(\dF(x), y)| =
        \begin{cases*}
            |x|, & if $x, y \in [-1, 0]$, \\
            y, & if $x, y \in [0, 1]$, \\
            |x + y|, & if $x \in [-1, 0]$ and $y \in [0, 1]$, \\
            0, & if $x \in [0,1]$ and $y \in [-1, 0]$.
        \end{cases*}
    \end{align*}
    Because the maximum of these values is~$1$,~$\codis(\dF, \dG) =
    1$.
    Lastly, because \eqref{eqn:max_in_sourcesink} holds for any pair of
    d-maps~$\dF, \dG$, and we found a specific pair for which the maximum
    equals~$1$,
    \begin{align*}
        \ddis(\X, \X^\ast)
        =\frac{1}{2} \inf_{\dF, \dG} \max \{\text{dis}(\dF), \text{dis}(\dG), \codis(\dF, \dG)\}
        = \frac{1}{2}.
    \end{align*}
\end{example}

\begin{proposition}
    Let $(\X,\dZ^X)$ and $(\Y,\dZ^Y)$ be bounded directed metric spaces. Then,
    \[
        \ddis(\X, \Y)
        \leq \frac{1}{2} \max \lbrace \diam(\X,\dZ^X), \diam(\Y,\dZ^Y) \rbrace.
    \]
\end{proposition}
\begin{proof}
    Consider $x_0$ in $X$ and $y_0$ in $Y$.
    Define the d-map $\dF_0\colon \X\to \Y$ as the constant map to $y_0$
    and $\dG_0\colon \Y\to \X$ as the constant d-map to $x_0$.
    We have that~$\dis(\dF_0)\le \diam(\X,\dZ^X)$,
    $\dis(\dG_0) \le \diam(\Y,\dZ^Y)$, and
    \begin{align*}
        \codis(\dF_0, \vec{G_0)}
        &=   \sup_{x\in\X, y\in Y}\lvert \dZ^{X}(x,x_0) - \dZ^{Y}(y,y_0)\rvert \\
        &\le \max \{\diam(\X,\dZ^X), \diam(\Y,\dZ^Y)\}.
    \end{align*}
    Thus, $\ddis(\X, \Y)\le \frac{1}{2} \max \{\diam(\X,\dZ^X), \diam(\Y,\dZ^Y)\}$.
\end{proof}

Consider two d-spaces $\X$ and $\Y$.
A \define{d-relation} $\mathcal{R}$
between $\X$ and $\Y$
is a relation between $X$ and $Y$ satisfying the following property:
for every
$(x,y),(x',y')$ in $\mathcal{R}$ and $\gamma_1\in \vec{P}(x,x')$,
there exists~$\gamma_2\in \vec{P}(y, y')$, and for every~$\gamma_2\in \vec{P}(y,
y')$, there
exists~$\gamma_1\in \vec{P}(x,x')$.
\begin{remark}
    The property characterizing d-relations is a reachability condition: If there is
    a d-path $\mu$ with $\mu(0)=a$ and $\mu(1)=b$, we say that $b$ is reachable
    from~$a$, or~$b$ is in the future of~$a$, or $a$ is in the past of $b$.
    This allows us to
    restate the property characterizing d-relations: For $(x,y)$ and $(x',y')$ in $\mathcal{R}$,
    $x'$ is reachable from~$x$ if and only if $y'$ is reachable from $y$ or
    equivalently $x'$ is in the future of~$x$ if and only if $y'$ is in the
    future of~$y$, and similarly for pasts.
\end{remark}

If the d-relation is also a correspondence, then we call it a
\define{d-correspondence}.

\begin{example}
    Consider the d-space $\vec{I}=(I, \vec{P}(I))$ with $\vec{P}(I)$ induced by
    the usual partial order.
    The relation
    $\mathcal{R}=\{(x,x)\in I\times I\}$ is a d-correspondence.
    The relation $\mathcal{R'}=\{(x,1-x)\in I\times I\}$, instead, is a
    correspondence, but not a d-correspondence.
\end{example}

The distortion is defined for d-correspondences with respect to the zigzag distance, that is
\[
\dis(\mathcal{R})=\sup_{(x,y),(x',y')\in \mathcal{R}}\lvert \dZ^X(x,x')-\dZ^Y(y,y')\rvert.
\]

\begin{definition}\label{d-corr-distance}
    The \define{d-correspondence distortion distance} between bounded directed metric
    spaces~$(\X, \dZ^X)$ and $(\Y, \dZ^Y)$ is defined as
     \begin{align}\label{def:d-corr-distance}
        \vec{d}_{\textup{c-dis}}(\X, \Y)
        = \frac{1}{2}\inf_{\mathcal{R}} \dis(\mathcal{R}),
    \end{align}
    where $\mathcal{R}$ varies over all the d-correspondences between $\X$
    and $\Y$.
\end{definition}

\begin{proposition}\label{prop_dis-distances-equal}
    For any
    directed metric spaces $\X$ and $\Y$,
    \begin{equation*}
        \dGH(\X,\Y) \leq \vec{d}_{\textup{c-dis}}(\X, \Y).
    \end{equation*}
\end{proposition}
\begin{proof}
    Assume $\vec{d}_{\textup{c-dis}}(\X, \Y) < \infty$, else the statement holds trivially.
    Because d-correspondences form a subset of correspondences, the result follows by \cref{thm:gh-vs-dgh} and \cref{eqn:kalton_bis}.
\end{proof}
We note that the relationship between $\vec{d}_{dis}(\vec{X}, \vec{Y})$ and $\vec{d}_{c-dis}(\vec{X}, \vec{Y})$ is still an open question. However, the next example shows that the two are not guaranteed to be equal.

\begin{remark}\label{rmk:cdis_greaterthan_dis}
    For $\X$ and $\Xrev$ as in \cref{ex:dK_greaterorequal_dGH},
    there exists no d-correspondence.
    Thus, $\vec{d}_{\textup{c-dis}}(\X, \Xrev)=\infty$; whereas,~$\ddis(\X,
    \Xrev)<\infty$.
\end{remark}

Recall that, for bounded metric spaces (not directed) $(X, d^X)$ and $(Y, d^Y)$,
we have the following equalities (\cref{eqn:kalton_bis} and \cref{eqn:kalton}):
\[
\GH (X, Y) =\frac{1}{2}\inf_{f, g} \max\{\dis(f), \dis(g), \codis(f, g)\} =\frac{1}{2}\inf_{\mathcal{R}} \dis(\mathcal{R}) .
\]
However, when considering bounded directed metric spaces $(\X, \dZ^X)$ and $(\Y, \dZ^Y)$
and d-isometries, the (co)distortion of d-maps, and d-correspondences between
them, these equalities do not necessarily hold. Instead, we get the following
inequalities:
\[
    \dGH(\X, \Y) \leq \ddis(\X, \Y),\;\;\;\; \dGH(\X, \Y) \leq \vec{d}_{\textup{c-dis}}(\X, \Y).
\]

\section{Discussion}
In this work, we introduced three alternative, non-equivalent, definitions of
Gromov--Hausdorff-type metrics between directed spaces. The question of which of
these definitions is most appropriate remains open and likely depends on the
specific application under consideration. We intend to explore this question
further in future work.

Our constructions are based on specific choices, leading to three notions of
directed Gromov--Hausdorff distance: one defined via \emph{d-isometries} between
directed metric spaces, and two based on \emph{distortion} of
\emph{d-maps} and \emph{d-correspondences}, respectively. These formulations
rely on particular assumptions regarding the continuity of maps that respect the
directed structures, as well as the choice of metrics used to define isometries,
distortions, and codistortions.  Future work will focus on investigating the
stability of the proposed directed distances and exploring further applications
that highlight their usefulness in applied settings. Additionally, we plan to
examine alternative definitions and analyze how they relate to each other and to
the classical Gromov--Hausdorff~distance.

\vspace{5mm}
\paragraph*{\textbf{Acknowledgements}}
This collaboration began at the WinCompTop3 workshop at EPFL in Lausanne,
Switzerland, which was partially funded by NSF CCF 2317401. The
authors thank the organizers---Erin Chambers, Heather Harrington,
Kathryn Hess, and Claudia
Landi---for organizing the productive workshop.  In addition, BTF was partially
supported by NSF CCF 2046730. LM was supported by the Mila EDI scholarship. FT was funded by the Knut and Alice Wallenberg Foundations and the WASP Postdoctoral Scholarship Program.
ZU is part of the Centre for TDA, supported by EPSRC grant EP/R018472/1.
We would like to thank anonymous reviewers whose careful reading have greatly improved this paper.

\bibliographystyle{plain}
\bibliography{references}

\end{document}